\documentclass{siamltex}
\usepackage{amsmath}
\usepackage{amssymb}
\usepackage{euscript}
\usepackage{xcolor}

\usepackage{ifpdf}
\ifpdf
  \usepackage[pdftex]{graphicx}
  \graphicspath{{Figures/}}
  \DeclareGraphicsExtensions{.pdf}
  \usepackage[pdftex]{hyperref}
\else
  \usepackage[dvips]{graphicx}
  \graphicspath{{Figures/}}
  \DeclareGraphicsExtensions{.eps}
  \usepackage[dvips]{hyperref}
\fi

\newtheorem{remark}{Remark}
\newtheorem{assumption}{Assumption}
\newtheorem{algorithm}{Algorithm}

\def\mathscr{\EuScript}

\newcommand{\defegal}{:=}                                   % D\'{e}finition

\newcommand{\bbR}{\mathbb{R}}                               % Nombres r\'{e}els

\newcommand{\norm}[1]{\left\|#1\right\|}                    % Norme
\newcommand{\projop}[1]{\mathrm{proj}_{#1}}                 % Op\'{e}rateur projection
\newcommand{\proj}[2]{\projop{#1}\left(#2\right)}           % Projection
\newcommand{\compo}{\mathop{\scriptstyle\circ}}             % Composition

\newcommand{\np}[1]{(#1)}                                   % Parenth\`{e}se normal
\newcommand{\bp}[1]{\big(#1\big)}                           % Parenth\`{e}se big
\newcommand{\Bp}[1]{\Big(#1\Big)}                           % Parenth\`{e}se Big
\newcommand{\bgp}[1]{\bigg(#1\bigg)}                        % Parenth\`{e}se bigg
\newcommand{\espace}[1]{\mathbb{#1}}                        % Espace de Hilbert
\newcommand{\fcara}[1]{\chi_{{}_{#1}}}                      % Fonction caract\'{e}ristique
\newcommand{\argmin}{\mathop{\arg\min}}                     % Arg-min
\newcommand{\gradi}[2][]{\nabla_{#1}#2}                     % Gradient partiel
\newcommand{\espacea}[1]{\mathbb{#1}}                       % Espace d'arriv\'{e}e
\newcommand{\espacef}[1]{\mathcal{#1}}                      % Espace fonctionnel
\newcommand{\tribu}[1]{\mathscr{#1}}                        % Tribu
\newcommand{\borel}[1]{\tribu{B}_{#1}^{\mathrm{o}}}         % Tribu borelienne
\newcommand{\omeg}{\Omega}                                  % espace du triplet
\newcommand{\trib}{\tribu{A}}                               % tribu  du triplet
\newcommand{\prbt}{\mathbb{P}}                              % proba  du triplet
\newcommand{\espe}{\mathbb{E}}                              % Symbole esp\'{e}rance
\makeatletter
\def\va@a{\boldsymbol{\va@arg^{\textstyle\text{\unboldmath$\scriptstyle\va@expo$}}_{\textstyle\text{\unboldmath$\scriptstyle\va@index$}}}}
\def\va#1{\def\va@expo{}\def\va@index{}\def\va@arg{\uppercase{#1}}%
  \@ifnextchar^{\va@h}{\@ifnextchar_\va@u\va@a}}
\def\va@h^#1{\def\va@expo{#1}\@ifnextchar_\va@hu\va@a}
\def\va@u_#1{\def\va@index{#1}\@ifnextchar^\va@uh\va@a}
\def\va@hu_#1{\def\va@index{#1}\va@a}
\def\va@uh^#1{\def\va@expo{#1}\va@a}
\makeatother

\newcommand{\normdelim}[1]{\np{#1}}                         % Taille ``normal''
\newcommand{\bigdelim}[1]{\bp{#1}}                          % Taille ``big''
\newcommand{\Bigdelim}[1]{\Bp{#1}}                          % Taille ``Big''
\newcommand{\biggdelim}[1]{\bgp{#1}}                        % Taille ``bigg''
\newcommand{\normdelims}[2]{\normdelim{#1\mid#2}}           % avec s\'{e}parateur
\newcommand{\bigdelims}[2]{\bigdelim{#1\ \big|\ #2}}        %
\newcommand{\Bigdelims}[2]{\Bigdelim{#1\ \Big|\ #2}}        %
\newcommand{\besp}[2][]{\espe_{#1}\bigdelim{#2}}            % Esp\'{e}rance big
\newcommand{\Besp}[2][]{\espe_{#1}\Bigdelim{#2}}            % Esp\'{e}rance Big
\newcommand{\bgesp}[2][]{\espe_{#1}\biggdelim{#2}}          % Esp\'{e}rance bigg
\newcommand{\nespc}[3][]{\espe_{#1}\normdelims{#2}{#3}}     % Esp\'{e}rance cond. normal
\newcommand{\bespc}[3][]{\espe_{#1}\bigdelims{#2}{#3}}      % Esp\'{e}rance cond. big
\newcommand{\Bespc}[3][]{\espe_{#1}\Bigdelims{#2}{#3}}      % Esp\'{e}rance cond. Big
\newcommand{\proscal}[2]{\left\langle#1\:,#2\right\rangle}  % Produit scalaire

\newcommand{\as}{\text{a.s.}}                               % ``a.s.''
\newcommand{\Pps}{\text{$\prbt$-}\as}                       % ``P-a.s.''

\newcommand{\grid}[1]{\boldsymbol{\lowercase{#1}}}            % Grille
\newcommand{\gridop}[1]{\mathfrak{#1}}                      % Operateur sur grille

\def\eqsepv{,\enspace}                                      % Virgule dans une \'{e}quation

\def\eqfinv{,}                                              % Virgule en fin d'\'{e}quation
\def\eqfinp{.}                                              % Point en fin d'\'{e}quation
\def\Ufe{\espacef{U}^{\mathrm{fe}}}
\def\Uas{\espacef{U}^{\mathrm{as}}}
\def\Ume{\espacef{U}^{\mathrm{me}}}

\def\Gammaas{\Gamma^{\mathrm{as}}}
\def\bGammaas{\va{\Gamma}^{\mathrm{as}}}

\def\Phias{\Phi^{\mathrm{as}}}

\def\Hcv{\espacef{H}^{\mathrm{cv}}}
\def\Hsp{\espacef{H}^{\mathrm{sp}}}
\def\Hfe{\espacef{H}^{\mathrm{fe}}}

\title{Particle Methods For Stochastic Optimal Control Problems}

\author{
Pierre Carpentier
\thanks{\'{E}cole Nationale Sup\'{e}rieure de Techniques Avanc\'{e}es,
32 boulevard Victor, 75739 Paris Cedex 15, France
({\tt Pierre.Carpentier@ensta.fr})}
\and Guy Cohen
\thanks{Universit\'{e} de Paris-Est, CERMICS, \'{E}cole des Ponts,
Champs sur Marne, 77455 Marne la Vall\'{e}e Cedex 2, France
({\tt guy.cohen@mail.enpc.fr})}
\and Anes Dallagi
\thanks{EDF R\&D,
1 avenue du G\'{e}n\'{e}ral de Gaulle, 92141 Clamart Cedex, France
({\tt anes.dallagi@edf.fr})}
}

\begin{document}

\maketitle

\begin{abstract}

When dealing with numerical solution of stochastic optimal control
problems, stochastic dynamic programming is the natural framework.
In order to try to overcome the so-called curse of dimensionality,
the stochastic programming school promoted another approach based
on scenario trees which can be seen as the combination of Monte
Carlo sampling ideas on the one hand, and of a heuristic technique
to handle causality (or nonanticipativeness) constraints on the
other hand.

However, if one considers that the solution of a stochastic
optimal control problem is a feedback law which relates control to
state variables, the numerical resolution of the optimization
problem over a scenario tree should be completed by a feedback
synthesis stage in which, at each time step of the scenario
tree, control values at nodes are plotted against corresponding
state values to provide a first discrete shape of this feedback
law from which a continuous function can be finally inferred. From
this point of view, the scenario tree approach faces an important
difficulty: at the first time stages (close to the tree root),
there are a few nodes (or Monte-Carlo particles), and therefore a
relatively scarce amount of information to guess a feedback law,
but this information is generally of a good quality (that is,
viewed as a set of control value estimates for some particular
state values, it has a small variance because the future of those
nodes is rich enough); on the contrary, at the final time stages
(near the tree leaves), the number of nodes increases but the
variance gets large because the future of each node gets poor (and
sometimes even deterministic).

After this dilemma has been confirmed by numerical experiments, we
have tried to derive new variational approaches. First of all, two
different formulations of the essential constraint of
nonanticipativeness are considered: one is called \emph{algebraic}
and the other one is called \emph{functional}. Next, in both
settings, we obtain optimality conditions for the corresponding
optimal control problem. For the numerical resolution of those
optimality conditions, an adaptive mesh discretization method is
used in the state space in order to provide information for
feedback synthesis. This mesh is naturally derived from a bunch of
sample noise trajectories which need not to be put into the form of
a tree prior to numerical resolution. In particular, an important
consequence of this discrepancy with the scenario tree approach is
that the same number of nodes (or points) are available from the
beginning to the end of the time horizon. And this will be
obtained without sacrifying the quality of the results (that is,
the variance of the estimates). Results of experiments with a
hydro-electric dam production management problem will be presented
and will demonstrate the claimed improvements.
\end{abstract}



\pagestyle{myheadings}
\thispagestyle{plain}
\markboth{A. DALLAGI AND G. COHEN  AND P. CARPENTIER}%
{PARTICLE METHODS FOR STOCHASTIC OPTIMIZATION PROBLEMS}

\section*{Introduction}

Taking into account uncertainties in the decision process has
become an important issue for all industries. Facing the market
volatilities, the weather whims and the changing policies and
regulatory constraints, decision makers have to find an \emph{optimal}
way to introduce their decision process into this uncertain framework.
One way to take uncertainties into account is to use the
\emph{stochastic optimization} framework. In this approach,
the decision maker makes his decision by
optimizing a mean value with regard to the multiple possible
scenarios weighted with a probability law.

Stochastic optimization problems often involve information
constraints: the decision maker makes his decision after getting
observations about the possible scenarios. In the literature, two
different communities have dealt with this information issue using
different modelling techniques, and therefore different solution
methods. The question is: knowing that stochastic optimization
problems are often infinite dimensional problems, how can we
implement tractable solution methods (which can be implemented
using a computer)? To answer this question each community brought
its own answers.

The stochastic programming community models the information
structure by scenarios trees: this involves Monte Carlo sampling
plus some manipulations of the sample trajectories to enforce the
tree structure. Then, tractable solution methods for stochastic
optimization problems consist in solving deterministic problems
over the decision trees.
We refer to \cite{Pflug01,HeitschRomisch03,DupacovaGroweKuskaRomisch03}
for further details on these methods. Nevertheless, stochastic
programming faces an important difficulty: due to the tree
structure, one has a few discretization points (nodes of the tree)
at the early stages, which could represent a serious handicap when
attempting to synthesizing a feedback law. On the contrary, when
approaching the last stages, one has a large number of
discretization points but with a future which may be almost or
completely deterministic.

The stochastic optimal control community uses special structures
of the stochastic optimization problems (time sequentiality and
state notions) to model the information constraints through a
functional interpretation that leads to the Dynamic Programming
Principle.
We refer to \cite{Bellman:1957,Bertsekas:1976,Bertsekas-Shreve:1996}
for further details on this method. However, stochastic dynamic
programming is also confronted with a serious obstacle known as the
``curse of dimensionality''. In fact, this method leads one to blindly
discretize the whole state space without taking into account a,
generally nonuniform, state distribution at the optimal solution.

This paper tries to bridge the gap between those two communities.
We propose a tractable solution method for
stochastic optimal control problems which makes use of the good
ideas of both Monte Carlo sampling and variational methods on the
one hand, functional handling of the information structure on the
other hand. Other related works
\cite{BrodieGlasserman97,ThenieVial06} are still closer to the
stochastic programming point of view. In our approach, the same
number of discretization (sample) points are used from the
beginning to the end of the time horizon, and this discretization
grid is adaptive: when transposed into the state space, it
reflects the optimal state distribution at each time stage.

The proposed method is of a variational nature in that it is based
on gradient calculations which involve the forward state variable
integration and the backward adjoint state (co-state) evaluation
as in the Pontryagin minimum principle. The purpose is to solve
the Kuhn-Tucker necessary optimality conditions of the considered
problem (including the information constraints).

This paper is organized as follows. In \S\ref{SecPrelim}, we discuss
how to model the information structure in stochastic optimization problems.
Then we derive optimality conditions for stochastic optimization problems
with information constraints. In \S\ref{OCalg}, these optimality
conditions are specialized to the situation of stochastic optimal control
problems. A special attention is paid to the so-called \emph{Markovian}
case in \S\ref{OCfunc}. The numerical implementation of a resolution
method based on Monte Carlo and functional approximations is presented
in \S\ref{adapdisc}. Finally, in \S\ref{SecTest}, a case study, namely
a hydro-electric dam management problem, illustrates the proposed
method, and compares it to the standard scenario tree approach.

\section{Preliminaries\label{SecPrelim}}

In this section, we present the main framework of this paper: how
to model stochastic optimization problems and how to represent
the information structure of such problems.
This preliminary section is directly inspired from the works
of the System and Optimization Working Group\footnote{SOWG,
\'{E}cole Nationale des Ponts et Chauss\'{e}es: Laetitia Andrieu, Kengy Barty,
Pierre Carpentier, Jean-Philippe Chancelier, Guy  Cohen, Anes Dallagi,
Michel De Lara, Pierre Girardeau, Babakar Seck, Cyrille Strugarek.}
\cite{TheseBarty,TheseStrugarek,TheseDallagi}.

In the sequel of the paper, the random variables, defined over
a probability space $(\omeg,\trib,\prbt)$, will be denoted using bold
letters (e.g. $\va{\xi} \in L^{2}(\omeg,\trib,\prbt;\Xi)$)
whereas their realizations will be denoted using normal letters
(e.g. $\xi \in \Xi$).

\subsection{Modelling information}

When unknown factors affect a process that we try to control,
we consider a set $\omeg$ of possible \emph{states of nature}
$\omega$ among which the true state $\omega_0$ is supposed to be.
Roughly speaking, an information structure is a partition
of $\omeg$ into a collection of subsets $G$. A posteriori
observations will help us to determine in which particular subset
of this partition the true state $\omega_0$ lies. Two extreme
cases may be mentioned.
\begin{enumerate}
\item If the partition is simply the crude partition
$\{\emptyset,\omeg\}$, then the a posteriori observations are
useless. We will refer to this situation as an \emph{open-loop}
information structure.
\item If the partition is the finest
possible one, namely all subsets $G$ are singletons, then the a
posteriori observations will tell us exactly what is the true
state of nature $\omega_0$. This is the situation of perfect
knowledge prior to making our decisions.
\end{enumerate}
In between, after the a posteriori observations become
available, we will remain with some uncertainty about the true $\omega_0$,
that is we will know in which particular $G$ the true state lies
but all $\omega$'s in that $G$ are still possible.
Moreover, in dynamic situations, new observations may become
available at each time stage and a partition of $\omeg$
corresponding to this information structure must be considered at
each time stage. We refer to the general case as a
\emph{closed-loop} decision process.

In order to have a framework in which various operations on
information structures become possible, probability theory
introduces so-called $\sigma$-algebras or $\sigma$-fields, random
variables and a pre-order relation between random variables. The
latter is called \emph{measurability}.\footnote{A random variable
$\va{y}_1$ is measurable with respect to another random variable
$\va{y}_2$ if and only if the $\sigma$-field generated by the first
is included into the one generated by the second:
$\sigma(\va{y}_1)\subset\sigma(\va{y}_2)$. In this paper, we will
not give further details on these definitions. Instead, we refer
the reader to \cite{Brei92} for further details on the probability
and measurability theory.}

We refer the reader to \cite{Rao04} for further details
on information structure and probability theory. As far as this
paper is concerned, we here state a result found in
\cite[Theorem~8~p.108]{Rao04}  giving the main properties of the
measurability relation between random variables.

\begin{proposition}[Measurability relation]
Let $(\omeg,\trib,\prbt)$ be a probability space,
let $\va{y}_1:\omeg\rightarrow\espacea{Y}_1$ and
$\va{y}_2:\omeg\rightarrow\espacea{Y}_2$ be two random variables
taking their values in $\espacea{Y}_1$ and $\espacea{Y}_2$ respectively.
The following statements are equivalent:
\begin{enumerate}
\item $\va{y}_1\preceq \va{y}_2$,
\item $\sigma(\va{y}_1)\subset\sigma(\va{y}_2)$,
\item $\exists \phi:\espacea{Y}_2\rightarrow \espacea{Y}_1$,
a measurable mapping such that $\va{y}_1=\phi\circ \va{y}_2$,
\item $\va{y}_1=\nespc{\va{y}_1}{\va{y}_2}\eqsepv\Pps$.
\end{enumerate}
\label{propmes}
\end{proposition}

\subsection{\label{MSOP}Modelling a stochastic optimization problem}

In this section, we consider two interpretations of a stochastic
optimization problem: an algebraic one in which the information
constraint is modelled through the standard measurability relation
(Statement~$1$ of Proposition~\ref{propmes}), and a functional
interpretation in which we use the functional equivalent of
measurability relation (Statement~$3$ of Proposition~\ref{propmes}).

\subsubsection{Algebraic interpretation}

Let $(\omeg,\trib,\prbt)$ be a probability space and let
$\va{\xi}:\omeg\rightarrow\Xi$ be a random variable taking
values in $\Xi:=\bbR^{d_\xi}$ (noise space). We denote by
$\espacea{U}:=\bbR^{d_u}$ the control space and by $\espacef{U}$
the functional space $L^2(\omeg,\trib,\prbt;\espacea{U})$.
The cost function $J:\espacef{U}\rightarrow\bbR$ is defined
as the expectation of a normal integrand\footnote{See
\cite[Definition~14.27]{Rockafellar-Wets:1998}.
We remind that the normal integrand assumption is done to ensure
measurability properties \cite[Proposition~14.28]{Rockafellar-Wets:1998}.}
$j:\espacea{U}\times\Xi\rightarrow\bbR$
\begin{equation*}
J(\va{u}) \defegal \besp{j(\va{u},\va{\xi})} =
\int_\omeg j\left(\va{u}(\omega),\va{\xi}(\omega)\right)\mathrm{d}\prbt(\omega)
\eqfinp
\end{equation*}
The feasible set $\Ufe\defegal\Uas\cap\Ume$ accounts for two different
 constraints:
\begin{itemize}
\item almost sure constraints:
\begin{equation}
\label{eq:as-const}
\Uas \defegal \left\{\va{u}\in\espacef{U} \eqsepv
\va{u}(\omega)\in\bGammaas(\omega) \eqsepv \Pps\right\} \eqfinv
\end{equation}
where $\bGammaas:\omeg\rightrightarrows\espacea{U}$ is a measurable
set-valued mapping (see \cite[Definition~14.1]{Rockafellar-Wets:1998})
which is convex and closed valued, and
\item measurability constraints:
\begin{equation}
\label{eq:mes-const}
\Ume\defegal\left\{\va{u}\in\espacef{U} \eqsepv
\va{u}\enspace\text{is}\enspace\tribu{G}\mathrm{-measurable}\right\} \eqfinv
\end{equation}
where $\tribu{G}$ is a given sub-$\sigma$-field of $\trib$.
\end{itemize}
From their respective definitions, it is easy to prove that $\Ume$ is a closed
subspace of $\espacef{U}$ (indeed $L^2(\omeg,\tribu{G},\prbt;\espacea{U})$)
and that $\Uas$ is a closed convex subset of $\espacef{U}$.

The optimization problem under consideration is to minimize the cost function
$J(\va{U})$ over the feasible subset $\Ufe$.
The first model representing the stochastic optimization problem is thus
\begin{equation}
\min_{\va{u}\in\espacef{U}} \;  J(\va{u})
\quad \text{s.t.} \quad \va{u} \in \Ufe \eqfinp
\label{GSP}
\end{equation}
%%\renewcommand{\theequation}{\arabic{equation}}
%%\addtocounter{equation}{-1}
\noindent
This interpretation is called algebraic as it uses an algebraic
relation (measurability) to define the information structure of
Problem~\eqref{GSP}.

\subsubsection{Functional interpretation}

According to Proposition~\ref{propmes}, a functional model for stochastic
optimization problems is available, in which the optimization is achieved
with respect to functions called \emph{feedbacks}. Indeed, let
\begin{itemize}
\item $\va{y}:\omeg\rightarrow\espacea{Y}$ be a random variable
(called \emph{the observation}) taking value in $\espacea{Y}:=\bbR^{d_y}$,
\item $\Phi$ be the space
$L^2(\espacea{Y},\borel{\espacea{Y}},\prbt_{\va{y}};\espacea{U})$
where $\borel{\espacea{Y}}$ is the Borel $\sigma$-field of $\espacea{Y}$
and $\prbt_{\va{y}}$ the image of the probability measure $\prbt$ by $\va{y}$,
\item $\Phias$ be the subset of $\Phi$ defined by
\begin{equation*}
\Phias \defegal \big\{\phi\in\Phi \eqsepv \phi\compo\va{y}(\omega)
\in \bGammaas(\omega) \eqsepv \Pps\big\} \eqfinp
\end{equation*}
\end{itemize}
We define the cost function $\bar{J} : \Phi \rightarrow \bbR$
as $\bar{J}(\phi) \defegal \Besp{j\big(\phi(\va{y}),\va{\xi}\big)}$.
We are interested in the functional optimization problem:
\begin{equation}
\min_{\phi\in\Phi}\bar{J} \; (\phi)
\quad \text{s.t.} \quad \phi \in \Phias \eqfinp
\label{FSP}
\end{equation}

\begin{proposition}
If $\sigma(\va{y})=\tribu{G}$, then Problem~\eqref{GSP} is equivalent
to Problem~\eqref{FSP} in the sense that if $\va{u}^\sharp$ is solution
of~\eqref{GSP} (resp. $\phi^\sharp$ solution of~\eqref{FSP}),
then there exists $\phi^\sharp$ solution of~\eqref{FSP} (resp.
$\va{u}^\sharp$ solution of~\eqref{GSP}) such that
$\va{u}^\sharp=\phi^\sharp(\va{y})\eqsepv\Pps$.
\label{propalgfunc}
\end{proposition}
\begin{proof}
This is a straightforward consequence of
Proposition~\ref{propmes} and of the definition of the feasible
sets in both problems. \qquad
\end{proof}

\subsection{\label{OCSOP}Optimality conditions for a stochastic optimization problem}

Previous works dealt with optimality conditions for stochastic
optimal control problems \cite{Romisch05,Hiriart82}.
This paragraph can be viewed as a slight extension of \cite{Romisch05}
when considering a stochastic optimization problem subject to almost-sure
\emph{and} measurability constraints.

We consider Problem~\eqref{GSP} presented in \S\ref{MSOP}.
We recall that the feasible set $\Ufe$ is the intersection
of a closed convex subset $\Uas$ and of a linear subspace $\Ume$.
In order to apply the results given in Appendix \ref{app-OC-Hilbert},
we first establish the following lemma.

\medskip

\begin{lemma}
Let $\Uas$ be the convex set defined by almost sure constraints
\eqref{eq:as-const} and let $\Ume$ be the linear subspace
defined by measurability constraints~\eqref{eq:mes-const}.
We assume that $\bGammaas$ is a $\tribu{G}$-mesurable and closed
convex valued mapping. Then
\begin{equation*}
\proj{\Uas}{\Ume} \subset \Ume \eqfinp
\end{equation*}
\label{lemprint}
\end{lemma}
\begin{proof}
We first prove that
$\big(\proj{\Uas}{\va{u}}\big)(\omega)=
 \proj{\bGammaas(\omega)}{\va{u}(\omega)} \eqsepv \Pps$. Indeed, let
$F(\va{V}) \defegal \frac{1}{2} \norm{\va{v}-\va{u}}_{\espacef{U}}^{2} +
 \fcara{\Uas}(\va{v})$. From the definition of almost sure constraints,
we have
 $F(\va{v}) = \int_{\omeg} f\big( \va{v}(\omega),\omega\big)
 \mathrm{d}\prbt(\omega)$, with
$f(v,\omega) \defegal \frac{1}{2} \norm{v-\va{u}(\omega)}_{\espacea{U}}^{2} +
 \fcara{\bGammaas(\omega)}(v)$.
By definition of the projection, $\proj{\Uas}{\va{u}}$ is solution
of the optimization problem
\begin{equation*}
\min_{\va{V}\in\espacef{U}}
\int_{\omeg} f\big( \va{v}(\omega),\omega\big)\mathrm{d}\prbt(\omega) \eqfinp
\end{equation*}
Using \cite[Theorem 14.60]{Rockafellar-Wets:1998} (interchange
of minimization and integration), we obtain
\begin{equation}
\big(\proj{\Uas}{\va{u}}\big)(\omega) \in
\argmin_{v\in\espacea{U}} f(v,\omega) =
\big\{ \proj{\bGammaas(\omega)}{\va{u}(\omega)} \big\} \eqsepv \Pps \eqfinv
\label{eq:prf-pra}
\end{equation}
hence the claimed property.

Let us consider any $\va{u}\in\Ume$.
We deduce from~\eqref{eq:prf-pra} that
\begin{equation*}
\big(\proj{\Uas}{\va{u}}\big)(\omega) =
\argmin_{v\in\bGammaas(\omega)} \frac{1}{2}
\norm{v-\va{u}(\omega)}_{\espacea{U}}^{2} \eqsepv \Pps \eqfinp
\end{equation*}
From \cite[Theorem 8.2.11]{AubinFrankowska90} (measurability of marginal
functions), we deduce from the $\tribu{G}$-measurability of both $\va{U}$
and $\bGammaas$ that the $\argmin$ function $\proj{\Uas}{\va{u}}$ is also
a $\tribu{G}$-measurable function, which means that
$\proj{\Uas}{\va{u}}\in\Ume$.\qquad
\end{proof}

\medskip

\noindent
The main result of this section is given by the following theorem.

\begin{theorem}
Assume that $\bGammaas$ is a $\tribu{G}$-mesurable and closed convex
valued mapping, that function $j$ is a normal integrand such that
$j(\cdot,\va{\xi})$ is differentiable $\Pps$ and that
$j'_u(\va{u}(\cdot),\va{\xi}(\cdot))\in\espacef{U}$
for all $\va{u}\in\espacef{U}$.
Let $\va{u}^{\sharp}$ be a solution of Problem~\eqref{GSP}. Then
\begin{equation}
\bespc{j'_u(\va{u}^{\sharp},\va{\xi})}{\tribu{G}} \in
-\partial\fcara{\Uas}(\va{u}^{\sharp}) \eqfinp
\label{eq:GOC}
\end{equation}
\label{theoCOpr}
\end{theorem}
\begin{proof}
The differentiability of $J:\espacef{U}\rightarrow\bbR$ is
a straightforward consequence of both the integral expression of $J$
and the differentiability assumption on $j$. Moreover the expression
of the derivative of $J$ is given by
\begin{equation*}
J'(\va{u})(\omega)=j'_u(\va{u}(\omega),\va{\xi}(\omega)) \eqsepv \Pps \eqfinp
\end{equation*}
Let $\va{u}^{\sharp}$ be a solution of~\eqref{GSP}.
From Lemma \ref{lemprint} and Proposition \ref{propCOpr} in the Appendix,
we obtain
\begin{equation*}
\proj{\Ume}{J'(\va{u}^{\sharp})}\in -\partial\fcara{\Uas}(\va{u}^{\sharp})
\eqfinv
\end{equation*}
this last expression being equivalent to~\eqref{eq:GOC} thanks to
the characterization of the conditional expectation as a projection
and the expression of $J'$. \qquad
\end{proof}

\medskip

Using the equivalent statements~\eqref{cond3}, the optimality conditions
given by Theorem~\ref{theoCOpr} can be reformulated as
\begin{equation*}
\va{u}^{\sharp} =
\proj{\Uas}{\va{u}^{\sharp}-\epsilon \proj{\Ume}{\nabla J(\va{u}^{\sharp})}}
\eqsepv \forall \epsilon > 0 \eqfinp
\end{equation*}
This last formulation can be used in practice to solve Problem
\eqref{GSP} using a projected gradient algorithm. The projection
over $\Uas$ involves random variables, but we saw in the proof
of Lemma~\ref{lemprint} that $\projop{\Uas}$ can be performed
in a pointwise manner ($\omega$ per $\omega$), so that
the implementation of the algorithm is effective.

\begin{remark}
The optimality conditions~\eqref{eq:GOC}, which are given in terms
of random variables, also have a pointwise interpretation.
As a matter of fact, using the equivalent statements~\eqref{cond123}
and the pointwise interpretation of $\projop{\Uas}$, it is straightforward
to prove that $\va{r}\in\partial\fcara{\Uas}(\va{u})$ implies
$\va{r}(\omega)\in\partial\fcara{\bGammaas(\omega)}(\va{u}(\omega))\eqsepv\Pps$.
\end{remark}

\section{\label{OCalg}Stochastic optimal control problems}

Stochastic optimal control problems rest upon the same framework as
stochastic optimization problems: they can be modelled as
closed-loop stochastic optimization problems with a sequential
time structure. In this section, we deal with the algebraic
interpretation of stochastic optimal control problems in a
discrete time framework. We will derive optimality conditions
following the same principle as in \S\ref{OCSOP}.

\subsection{Problem formulation}

We consider a stochastic optimal control problem in discrete time,
$T$ denoting the time horizon. At each stage $t=0,\ldots,T$,
we denote by $\espacea{W}_t:=\bbR^{d_{w_t}}$ the noise space at time $t$.
Let $\espacef{W}_t$ be a set of random variables defined on
$(\omeg,\trib,\prbt)$ and taking values in $\espacea{W}_t$, and let
$\espacef{W}\defegal\espacef{W}_0\times\cdots\times\espacef{W}_T$.
The noise process of the problem is a random vector $\va{w}\in\espacef{W}$
such that $\va{w}=(\va{w}_0,\ldots,\va{w}_T)$, with $\va{w}_t\in \espacef{W}_t$.

At each stage $t=0,\ldots,T-1$, we denote by
$\espacea{U}_t \defegal \bbR^{d_{u_t}}$ the control space and by
$\espacef{U}_t \defegal L^2(\omeg,\trib,\prbt;\espacea{U}_t)$
the space of square integrable random variables taking values in
$\espacea{U}_t$. At $t$, the decision maker makes a decision (a control)
$\va{u}_t\in\espacef{U}_t$.
Let $\espacef{U} \defegal \espacef{U}_0\times\cdots\times\espacef{U}_{T-1}$.
The control process of the problem is a random vector $\va{u}\in\espacef{U}$
such that $\va{u}=(\va{u}_0,\ldots,\va{u}_{T-1})$,
with $\va{u}_t\in \espacef{U}_t$.
We also assume that each control variable $\va{u}_{t}$ is subject to
almost-sure constraints. More precisely, for all $t=0,\ldots,T-1$,
let $\bGammaas_t: \Omega\rightrightarrows\espacea{U}_t$ be a set-valued
mapping (random set), let $\Uas_t\defegal\left\{\va{u}_t\in\espacef{U}_t,\;
 \va{u}_t(\omega)\in \bGammaas_t(\omega) \eqsepv \Pps\right\}$ and let
$\Uas\defegal\Uas_0\times\cdots\times\Uas_{T-1}$.
The almost-sure  constraints writes
\begin{equation}
\va{u}_t\in \Uas_t \eqsepv \forall t=0,\ldots,T-1 \eqfinp
\label{ascons}
\end{equation}

\subsubsection{Dynamics}

At each stage $t=0,\ldots,T$, we denote by $\espacea{X}_t:=\bbR^{d_{x_t}}$
the state space at time $t$. Let $\espacef{X}_t$ be a set of random variables
defined on $(\omeg,\trib,\prbt)$ and taking values in $\espacea{X}_t$,
and let $\espacef{X}\defegal\espacef{X}_0 \times\ldots\times\espacef{X}_T$.
The state process of the problem is a random vector $\va{x}\in\espacef{X}$
such that $\va{x}=(\va{x}_0,\ldots,\va{x}_T)$, with $\va{x}_t\in \espacef{X}_t$.
It arises from the dynamics $f_t:\espacea{X}_t\times
 \espacea{U}_t\times\espacea{W}_{t+1}\rightarrow\espacea{X}_{t+1}$
of the system, namely
\begin{subequations}
\begin{align}
\va{x}_0    & =\va{w}_0 \eqsepv \Pps \eqfinv\\
\va{x}_{t+1}& =f_t(\va{x}_t,\va{u}_t,\va{w}_{t+1}) \eqsepv \Pps
  \quad \forall t=0,\ldots, T-1 \eqfinp
\end{align}
\label{SOCDyn}
\end{subequations}

\begin{remark}
The shift on the time index between the random variable $\va{x}_t$
and the control $\va{u}_t$ on the one hand and the noise $\va{w}_{t+1}$
on the other hand enlightens the fact that we are in the so-called
\emph{decision-hazard} scheme: at time $t$, the decision maker chooses
a control variable $\va{u}_{t}$ before having any information about
the noise $\va{w}_{t+1}$ which affects the dynamics of the system.
\end{remark}

\subsubsection{Cost function}

We consider at each stage $t=0,\ldots,T-1$, an ``integral'' cost function
$L_t:\espacea{X}_t\times\espacea{U}_t\times\espacea{W}_{t+1}\rightarrow\bbR$.
Moreover, at the final stage $T$, we consider a final cost function
$K:\espacea{X}_T\rightarrow\bbR$. Summing up all these instantaneous costs,
we obtain the overall cost function of the problem
\begin{equation*}
\widetilde{\jmath}(x,u,w)\defegal\sum_{t=0}^{T-1}L_t(x_t,u_t,w_{t+1})+K(x_T)
\eqfinv
\end{equation*}
and the decision maker has to choose the control variables $\va{u}_t$ in order
to minimize the overall cost expectation
\begin{equation}
\widetilde{J}(\va{x},\va{u})\defegal
\bgesp{\sum_{t=0}^{T-1}L_t(\va{x}_t,\va{u}_t,\va{w}_{t+1})+K(\va{x}_T)}
\eqfinp
\label{SOCCost}
\end{equation}

\subsubsection{Information structure}

Generally speaking, the information available on the system at time $t$
is modelled as a random variable $\va{y}_t:\Omega\rightarrow\espacea{Y}_t$,
where $\espacea{Y}_t:=\bbR^{d_{y_t}}$ is the observation space.
Let $\espacef{Y}_t$ be the set of random variables taking their values in
$\espacea{Y}_t$. We suppose that there exists an observation mapping
$\widetilde{h}_t:\espacea{W}_0\times\cdots\times\espacea{W}_T\rightarrow
 \espacea{Y}_t$ such that
$\va{y}_t=\widetilde{h}_t(\va{w}_0,\ldots,\va{w}_T) \eqsepv \Pps$.
For all $t=0,\ldots,T$, we denote by $\tribu{G}_t=\sigma(\va{y}_t)$
the sub-$\sigma$-field of $\trib$ generated by the random variable
$\va{y}_t$:\footnote{Note that the $\sigma$-fields $\tribu{G}_t$'s
are ``fixed'', in the sense that they do not depend on the control
variable $\va{u}$.}
\begin{equation*}
\tribu{G}_t = \sigma \big( \widetilde{h}_t(\va{w}_0,\ldots,\va{w}_T) \big)
\eqfinp
\end{equation*}
The decision maker knows $\va{y}_t$ when choosing
the appropriate control $\va{u}_t$ at time $t$, so that the information
constraint is $\va{u}_t \preceq \va{y}_t$.
Using the notations $\Ume_t=\left\{\va{u}_t\in\espacef{U}_t \eqsepv
\va{u}_t\;\;\text{is}\;\;\tribu{G}_t\mathrm{-measurable}\right\}$
and $\Ume=\Ume_0\times\cdots\times\Ume_{T-1}$, the measurability
constraints of the problem writes
\begin{equation}
\va{u}_t \in \Ume_t \eqsepv \forall t=0,\ldots,T-1 \eqfinp
\label{SOCInfo}
\end{equation}

We also introduce the sequence of $\sigma$-fields
$\big(\tribu{F}_t\big)_{t=0,\ldots,T}$ associated with the noise process
$\va{w}$, where $\tribu{F}_t$ is the $\sigma$-field generated by the noises
prior to $t$:
\begin{equation*}
\tribu{F}_t=\sigma\big(\va{w}_0,\ldots,\va{w}_t\big) \eqsepv \forall t=0,\ldots,T
\eqfinp
\end{equation*}
This sequence is a filtration as it satisfies the inclusions
$\tribu{F}_0\subset\tribu{F}_1\subset\ldots\subset\tribu{F}_T\subset\trib$.
When $\tribu{G}_t=\tribu{F}_t$, we are in the case of complete causal
information.

\subsubsection{Optimization problem}

According to the notation given in the previous paragraphs,
the stochastic optimal control problem we consider here is
\begin{subequations}
\begin{equation}
\min_{(\va{u},\va{x})\in\espacef{U}\times\espacef{X}} \quad
\bgesp{\sum_{t=0}^{T-1}L_{t}(\va{x}_{t},\va{u}_{t},\va{w}_{t+1})+K(\va{x}_{T})}
\eqfinv
\end{equation}
subject to both the dynamics constraints~\eqref{SOCDyn}
and the control constraints
\begin{equation}
\va{u}_t\in \Uas_t \cap \Ume_t \eqsepv \forall t=0,\ldots,T-1
\eqfinp
\end{equation}
\label{SOC}
\end{subequations}
We follow here the algebraic interpretation of an optimization problem
given in \S\ref{MSOP}, as the information structure is defined using
a measurability relation.

In the way we modelled Problem~\eqref{SOC}, we have to minimize the cost
function $\widetilde{J}(\va{u},\va{x})$ with respect to both $\va{u}$
and $\va{x}$ under the dynamics constraints~\eqref{SOCDyn}.
But the state variables $\va{x}_{t}$ are in fact intermediary variables
and it is possible to eliminate them by recursively incorporating the
dynamics equations~\eqref{SOCDyn} into the cost function~\eqref{SOCCost}.
The resulting cost function $J$ only depends on the control variables
$\va{u}_{t}$. Its expression is $J(\va{u}) = \besp{j(\va{u},\va{W})}$,
with
\begin{equation*}
j(u,w) =
\sum_{t=0}^{T-1}\widetilde{L}_t(u_0,\ldots,u_t,w_0,\ldots,w_{t+1})+
                \widetilde{K}(u_0,\ldots,u_{T-1},w_0,\ldots,w_T) \eqfinp
\end{equation*}
In this setting, Problem~\eqref{SOC} is equivalent to
\begin{equation}
\min_{\va{u}\in\espacef{U}} \quad J(\va{u})
\quad \text{s.t.} \quad \va{u}\in\Uas\cap\Ume \eqfinp
\label{SO}
\end{equation}

The derivatives of the cost function $J$ can be obtained from
the derivatives of $\widetilde{J}$ using the well known adjoint
state method. This is based on the following result stated here
without proof.

\begin{proposition}
Assuming that functions $f_{t}$ and $L_{t}$ are continuously
differentiable with respect to their first two arguments and
that function $K$ is continuously differentiable, the partial
derivatives of $j$ with respect to $u$ are given by
\begin{equation*}
(j)'_{u_t}(u,w) =
(L_t)_{u_t}'(x_t,u_t,w_{t+1})+\lambda_{t+1}^\top(f_t)_{u_t}'(x_t,u_t,w_{t+1})
\eqfinv
\end{equation*}
where the state vector $(x_0,\ldots,x_T)$ satisfies
the forward dynamics equation
\begin{equation*}
x_{0}=w_{0} \eqsepv x_{t+1}=f_t(x_t,u_t,w_{t+1}) \eqfinv
\end{equation*}
whereas the adjoint state (or co-state) vector $(\lambda_0,\ldots,\lambda_T)$
is chosen to satisfy the backward dynamics equation
\begin{equation*}
\lambda_T=K'^\top(x_T) \eqsepv
\lambda_t=(L_t)'^\top_{x}(x_t,u_t,w_{t+1})+(f_t)'^\top_{x}(x_t,u_t,w_{t+1})
\lambda_{t+1} \eqfinp
\end{equation*}
\label{propcuisine}
\end{proposition}

\subsubsection{Assumptions}

In order to derive optimality conditions for Problem~\eqref{SOC},
we make the following assumptions.

\vspace{0.2cm}

\begin{assumption}[Constraints structure]~
\begin{enumerate}
\item[] $\bGammaas_t$ are closed convex set-valued mappings,
      $\forall t=0,\ldots,T-1$.
\end{enumerate}
\label{AssStruct}
\end{assumption}

\vspace{0.2cm}

\pagebreak

\begin{assumption}[Differentiability]
\begin{enumerate}
\item Functions $f_{t}$ (dynamics) and $L_{t}$ (cost) are continuous
      differentiable with respect to their first two arguments
      (state and control), $\forall t=0,\ldots,T-1$.
\item Function $K$ (final cost) is continuously differentiable.
\item Functions $L_t$ and $f_t$ are normal integrands,
      $\forall t=0,\ldots,T-1$.
\item The derivatives of $f_t$, $L_t$ and $K$ are square integrable,
      $\forall t=0,\ldots,T-1$.
\end{enumerate}
\label{AssDiff}
\end{assumption}

\vspace{0.2cm}

\begin{assumption}[Nonanticipativity and measurability]~
\begin{enumerate}
\item $\tribu{G}_t\subset\tribu{F}_t \eqsepv \forall t=0,\ldots,T$.
\item Mappings $\bGammaas_t$ are $\tribu{G}_t$-measurable,
      $\forall t=0,\ldots,T-1$.
\end{enumerate}
\label{AssNonAnt}
\end{assumption}

\vspace{0.2cm}

Assumption~\ref{AssDiff} will allow us all integration and derivation
operations needed for obtaining optimality conditions.
The measurability condition on $\bGammaas_t$ in Assumption~\ref{AssNonAnt}
expresses that the almost-sure constraints must at least have the same
measurability as the decision variables.
The first condition in Assumption~\ref{AssNonAnt} expresses
the causality of the problem: the decision maker has no access to information
in the future! Under this assumption, there exists a measurable mapping
$h_t:\espacea{W}_0\times\cdots\times\espacea{W}_t\rightarrow\espacea{Y}_t$
such that the information variable $\va{y}_t$ writes
\begin{equation*}
\va{y}_t=h_t(\va{w}_0,\ldots,\va{w}_t) \eqsepv \Pps \eqfinp
\end{equation*}

From Assumption~\ref{AssStruct} and using Lemma~\ref{lemprint}
the following property is readily available.
\begin{proposition}[Constraints structure]
For all $t=0,\ldots,T-1$,
\begin{enumerate}
\item $\Uas_t$ is a closed convex subset of $\espacef{U}_t$,
\item $\proj{\Uas_t}{\Ume_t}\subset\Ume_t$.
\end{enumerate}
\label{PrpStruct}
\end{proposition}

\subsection{Optimality conditions in stochastic optimal control problems}

We present here necessary optimality conditions for the stochastic optimal
control problem~\eqref{SOC}, which are an extension of the conditions
given in Theorem~\ref{theoCOpr}.

\subsubsection{Non-adapted optimality conditions}

A first set of optimality conditions is given in the next theorem.

\begin{theorem}
Let the two random processes $(\va{x}_t)_{t=0,\ldots,T}$ $\in\espacef{X}$
and $(\va{u}_t)_{t=0,\ldots,T-1}$ $\in\espacef{U}$ be a solution of
Problem~\eqref{SOC}. Suppose that Assumption~\ref{AssStruct},
\ref{AssDiff} and \ref{AssNonAnt} are satisfied. Then, there exists
a random process $(\va{\lambda}_t)_{t=0,\ldots,T}\in\espacef{X}$ such that,
for all $t=0,...,T-1$,
\begin{subequations}
\begin{align}
& \va{x}_0=\va{w}_0
  \eqfinv \label{eqPMS1-a} \\
& \va{x}_{t+1} = f_{t}(\va{x}_{t},\va{u}_{t},\va{w}_{t+1})
  \eqfinv \label{eqPMS1-b} \\
& \nonumber \\
& \va{\lambda}_{T} = K'^\top(\va{x}_{T})
  \eqfinv \label{eqPMS1-c} \\
& \va{\lambda}_{t} = (L_{t})'^\top_{x}(\va{x}_{t},\va{u}_{t},\va{w}_{t+1})+
                     (f_{t})'^\top_{x}(\va{x}_{t},\va{u}_{t},\va{w}_{t+1})
                     \va{\lambda}_{t+1}
  \eqfinv \label{eqPMS1-d} \\
& \nonumber \\
& \bespc{(L_{t})'_{u}(\va{x}_{t},\va{u}_{t},\va{w}_{t+1})+
          \va{\lambda}_{t+1}^{\top}
          (f_{t})'_{u}(\va{x}_{t},\va{u}_{t},\va{w}_{t+1})}
        {\tribu{G}_{t}} \in - \partial\fcara{\Uas_t}(\va{u}_t)
  \eqfinp \label{eqPMS1-e}
\end{align}
\label{eqPMS1}
\end{subequations}
\label{thmPMS1}
\end{theorem}

\begin{proof}
From the equivalence between Problem~\eqref{SOC} and Problem~\eqref{SO},
we obtain, using Theorem~\ref{theoCOpr}, that the solution $\va{u}$
satisfies
\begin{equation*}
\proj{\Ume}{J'(\va{u})}\in -\partial\fcara{\Uas}(\va{u}) \eqfinv
\end{equation*}
and therefore
$\bespc{j'_{u_t}(\va{u},\va{w})}{\tribu{G}_t}\in-\partial\fcara{\Uas_t}(\va{u}_t)$
for all $t=0,\ldots,T-1$.
The desired result follows from Proposition~\ref{propcuisine}. \qquad
\end{proof}

The  conditions given by Theorem~\ref{thmPMS1} are called
\emph{non-adapted optimality conditions} because the dual random process
$(\va{\lambda}_t)_{t=0,\ldots,T}$ is not adapted to the natural filtration
$(\tribu{F}_t)_{t=0,\ldots,T}$, that is, $\va{\lambda}_t$ generally depends
on the future. We will see in the next section that similar optimality
conditions can be written with help of an adapted dual random process.

\subsubsection{Adapted optimality conditions}

The following theorem presents optimality conditions involving
an adapted dual random process.

\begin{theorem}
Let the two random processes $(\va{x}_t)_{t=0,\ldots,T}$ $\in\espacef{X}$
and $(\va{u}_t)_{t=0,\ldots,T-1}$ $\in\espacef{U}$ be a solution of
Problem~\eqref{SOC}. Assume that Assumption~\ref{AssStruct}, \ref{AssDiff}
and \ref{AssNonAnt} are satisfied. Then, there exists a process
$(\va{\Lambda}_t)_{t=0,\ldots,T}\in\espacef{X}$ adapted to the filtration
$(\tribu{F}_{t})_{t=0,\ldots,T} $ such that, for all $t=0,...,T-1$,
\begin{subequations}
\begin{align}
& \va{x}_0=\va{w}_0 \eqfinv \\
& \va{x}_{t+1} = f_{t}(\va{x}_{t},\va{u}_{t},\va{w}_{t+1}) \eqfinv \\
& \nonumber \\
& \va{\Lambda}_{T} = K'^\top(\va{x}_{T}) \eqfinv \\
& \va{\Lambda}_{t} =
  \bespc{(L_{t})'^\top_{x}(\va{x}_{t},\va{u}_{t},\va{w}_{t+1})+
         (f_{t})'^\top_{x}(\va{x}_{t},\va{u}_{t},\va{w}_{t+1})\va{\Lambda}_{t+1}}
        {\tribu{F}_t} \eqfinv \\
& \nonumber \\
& \bespc{(L_{t})'_{u}(\va{x}_{t},\va{u}_{t},\va{w}_{t+1})+
          \va{\Lambda}_{t+1}^\top(f_{t})'_{u}(\va{x}_{t},\va{u}_{t},\va{w}_{t+1})}
        {\tribu{G}_{t}} \in-\partial\fcara{\Uas_t}(\va{u}_t) \eqfinp
\end{align}
\label{eqPMS2}
\end{subequations}
\label{thmPMS2}
\end{theorem}

\begin{proof}
All assumptions of Theorem \ref{thmPMS1} are met, so that there exists
a random process $(\va{\lambda}_{t})_{t=0,\ldots,T}$ satisfying~\eqref{eqPMS1}.
Define for all $t$ the random variable $\Lambda_{t}$ by
\begin{equation*}
\va{\Lambda}_{t} \defegal \bespc{\va{\lambda}_{t}}{\tribu{F}_{t}} \eqfinp
\end{equation*}
By construction, the process $(\va{\Lambda}_{t})_{t=0,\ldots,T}$
is adapted to the filtration $(\tribu{F}_{t})_{t=0,\ldots,T}$.
At stage $T$, we have
\begin{equation*}
\va{\Lambda}_{T} = \nespc{\va{\lambda}_T}{\tribu{F}_T}
                 = \bespc{K'^\top(\va{x}_T)}{\tribu{F}_T}
                 = K'^\top(\va{x}_T) \eqfinv
\end{equation*}
because $\va{x}_T$ is $\tribu{F}_T$-measurable.
For all $t=T-1,\ldots,0$, using the law of total expectation
$\nespc{\cdot}{\tribu{F}_t}=\bespc{\nespc{\cdot}{\tribu{F}_{t+1}}}{\tribu{F}_t}$
and since all variables $\va{x}_t$, $\va{u}_t$ and $\va{w}_{t+1}$
are $\tribu{F}_{t+1}$-measurable, we deduce from~\eqref{eqPMS1} that
\begin{align*}
\va{\Lambda}_t
& = \bespc{(L_t)'^\top_x(\va{x}_t,\va{u}_t,\va{w}_{t+1})+
           (f_t)'^\top_x(\va{x}_t,\va{u}_t,\va{w}_{t+1})
           \nespc{\va{\lambda}_{t+1}}{\tribu{F}_{t+1}}}
          {\tribu{F}_t} \eqfinv \\
& = \bespc{(L_t)'^\top_x(\va{x}_t,\va{u}_t,\va{w}_{t+1})+
           (f_t)'^\top_x(\va{x}_t,\va{u}_t,\va{w}_{t+1})
           \va{\Lambda}_{t+1}}
          {\tribu{F}_t} \eqfinv
\end{align*}
hence the adapted backward dynamics equations given in~\eqref{eqPMS2}.

Assumption~\ref{AssNonAnt} implies that
$\tribu{G}_t\subset\tribu{F}_t\subset\tribu{F}_{t+1}$.
Using $\nespc{\cdot}{\tribu{G}_t}=\bespc{\nespc{\cdot}{\tribu{F}_{t+1}}}{\tribu{G}_t}$,
and the measurability properties of $\va{x}_t$, $\va{u}_t$ and $\va{w}_{t+1}$,
the last optimality condition in~\eqref{eqPMS1} becomes
\begin{equation*}
\bespc{(L_t)'_u(\va{x}_t,\va{u}_t,\va{w}_{t+1})+
       \nespc{\va{\lambda}_{t+1}^\top}{\tribu{F}_{t+1}}
       (f_t)'_u(\va{x}_t,\va{u}_t,\va{w}_{t+1})}
      {\tribu{G}_t} \in -\partial\fcara{\Uas_t}(\va{u}_t) \eqfinv
\end{equation*}
hence the last optimality condition given in~\eqref{eqPMS2}. \qquad
\end{proof}

Note that in the optimality conditions given in Theorem~\ref{thmPMS2},
at each stage $t$, the gradient is projected over the subspace generated
by the observation $\sigma$-field $\tribu{G}_t$, whereas the adapted dual
random variable is projected over the subspace generated by $\tribu{F}_t$
which corresponds to the natural filtration of Problem~\eqref{SOC}.

\subsection{\label{KKTX}Optimality conditions in the Markovian case}

As far as information structure is concerned, the optimality
conditions~\eqref{eqPMS1} and~\eqref{eqPMS2} were obtained assuming
only nonanticipativeness (see Assumption~\ref{AssNonAnt}) and the fact
that the observation $\sigma$-fields $(\tribu{G}_t)_{t=0,\dots,T}$ do
not depend on the decisions $(\va{U}_t)_{t=0,\dots,T-1}$. The last
assumption allows us to avoid the so-called ``dual effect of control''
(see \cite{barty03} for further details). Another feature often
available in practice for the information structure is the
``perfect memory'' property, which intuitively means that
the information is not lost over time. The last property implies
that $(\tribu{G}_t)_{t=0,\dots,T}$ is a filtration, namely
\begin{equation*}
\tribu{G}_t\subset\tribu{G}_{t+1}\eqsepv \forall t=0,\dots T-1 \eqfinp
\end{equation*}
We will assume in the sequel that the perfect memory property holds,
and we will moreover assume complete causal noise observation, so that
$\tribu{G}_t=\tribu{F}_t \eqsepv \forall t=0,\ldots,T$.\footnote{Note
that complete causal noise observation implies the perfect memory
property, as far as $(\tribu{F}_t)_{t=0,\dots,T}$ is a filtration.}
Then both optimality conditions~\eqref{eqPMS1} and~\eqref{eqPMS2}
involve conditional expectations with respect to a random observation
variable the dimension of which increases with time (a new noise
random variable becomes available at each stage $t$).
This leads to a computational difficulty, of the same nature as
the so-called curse of dimensionality. To address this difficulty,
it would be an easier situation to have a constant dimension for
the observation space, as in the stochastic dynamic programming
principle \cite{Bertsekas:1976,Bertsekas-Shreve:1996} when the optimal
control at $t$ only depends on the state variable at the same time stage.
We thus consider new (more restrictive) assumptions which match the
stochastic optimal control framework an lead us to the desired situation.

\vspace{0.2cm}

\begin{assumption}[Markovian case]
\begin{enumerate}
\item $\tribu{G}_t=\tribu{F}_t \eqsepv \forall t=0,\ldots,T$
      (perfect memory and causal noise observation).
\item The random variables $\va{w}_{0},\ldots,\va{w}_{T}$ are independent
      (white noise).
\item The mappings $\bGammaas_t$ are constant (deterministic constraints):\\
      $\forall t=0,\ldots,T-1 \eqsepv \exists \: \Gammaas_t\subset\espacea{U}_t
       \eqsepv \bGammaas_t(\omega) = \Gammaas_t \eqsepv \Pps$.
\end{enumerate}
\label{AssMark}
\end{assumption}

\vspace{0.2cm}

The standard formulation of a stochastic optimal control
problem in the Markovian case is to assume that the state is
completely and perfectly observed. The problem formulation is accordingly
\begin{subequations}
\begin{equation}
\min_{(\va{u},\va{x})\in\espacef{U}\times\espacef{X}} \quad
\bgesp{\sum_{t=0}^{T-1}L_{t}(\va{x}_{t},\va{u}_{t},\va{w}_{t+1})+K(\va{x}_{T})}
\eqfinv
\end{equation}
subject to both the dynamics constraints~\eqref{SOCDyn}
and the control constraints
\begin{equation}
\va{u}_t \in \Gammaas_t \eqsepv \Pps \quad \text{and} \quad
\va{u}_t \preceq \va{x}_t \eqsepv \forall t=0,\ldots,T-1
\eqfinp
\end{equation}
\label{SOCX}
\end{subequations}

We now consider the optimality conditions~\eqref{eqPMS1} and~\eqref{eqPMS2}
and we specialize them to the Markovian case.

\subsubsection{Markovian case: non-adapted optimality conditions}

We present a non-adapted version of the optimality
conditions of Problem~\eqref{SOC} with Markovian assumptions. We
begin by presenting a result inspired by the stochastic dynamic
programming principle.

\begin{theorem}
Suppose that Assumptions~\ref{AssStruct}, \ref{AssDiff} and
\ref{AssMark} are fulfilled, and assume that there exist
two random processes $(\va{u}_t)_{t=0,\ldots,T-1}\in\espacef{U}$
and $(\va{x}_t)_{t=0,\ldots,T}\in\espacef{X}$ solution of Problem
\eqref{SOC}. Then there exists a process
$(\va{\lambda}_t)_{t=0,\ldots,T-1}\in\espacef{X}$ satisfying
\eqref{eqPMS1} and such that, for all $t=0,\ldots,T-1$,
\begin{align*}
& (a) \enspace
  \va{\lambda}_{t+1}\preceq
  \big(\va{x}_{t+1},\va{w}_{t+2},\ldots,\va{w}_T\big) \eqsepv \\
& (b) \enspace
  \va{u}_t\preceq \va{x}_t \eqsepv \\
& (c) \enspace
  \bespc{(L_t)'_u(\va{x}_t,\va{u}_t,\va{w}_{t+1})+
          \va{\lambda}_{t+1}^\top(f_t)'_u(\va{x}_t,\va{u}_t,\va{w}_{t+1})}
        {\tribu{F}_t} = \\
& \quad\quad\quad\quad\quad
  \bespc{(L_t)'_u(\va{x}_t,\va{u}_t,\va{w}_{t+1})+
         \va{\lambda}_{t+1}^\top(f_t)'_u(\va{x}_t,\va{u}_t,\va{w}_{t+1})}
        {\va{x}_t} \eqsepv \Pps \eqfinp
\end{align*}
\label{thmPMS1X}
\end{theorem}

\begin{proof}
First, Assumption~\ref{AssNonAnt} being implied by  Assumption~\ref{AssMark},
the existence of the process $(\va{\lambda}_t)_{t=0,\ldots,T}\in\espacef{X}$
satisfying \eqref{eqPMS1} is given by Theorem~\ref{thmPMS1}.
Denoting by $H_{t}$ the Hamiltonian at time $t$, namely
$H_{t}(x,u,w,\lambda) = L_{t}(x,u,w) + \lambda^{\top}f_{t}(x,u,w)$,
optimality conditions \eqref{eqPMS1-d} and \eqref{eqPMS1-e} write
\begin{subequations}
\begin{align}
& \va{\lambda}_{t} = (H_{t})'^\top_{x}
                     (\va{x}_{t},\va{u}_{t},\va{w}_{t+1},\va{\lambda}_{t+1})
  \eqfinv \label{eqHPMS1-d} \\
& \bespc{(H_{t})'_{u}(\va{x}_{t},\va{u}_{t},\va{w}_{t+1},\va{\lambda}_{t+1})}
        {\tribu{F}_{t}} \in-\partial\fcara{\Uas_t}(\va{u}_t)
        \eqfinp \label{eqHPMS1-e}
\end{align}

The proof of statements $(a)$ and $(b)$ is obtained by induction.
For the sake of simplicity, we first prove the result when
$\Uas_t=\espacef{U}$, so that \eqref{eqHPMS1-e} reduces to the
equality condition:
\begin{equation}
\label{eqHPMS1-eq1}
\bespc{(H_{t})'_{u}(\va{x}_{t},\va{u}_{t},\va{w}_{t+1},\va{\lambda}_{t+1})}
      {\tribu{F}_{t}} = 0 \eqfinp
\end{equation}
\end{subequations}

\begin{itemize}
\item At stage $T$, we know from \eqref{eqPMS1-c} that
$\va{\lambda}_T = \mu_{T-1}(\va{x}_T)$, with $\mu_{T-1}$ being
a measurable function,\footnote{In fact, $\mu_{T-1}={K'}^{\top}$.
From Assumption \ref{AssDiff}, $\mu_{T-1}$ is a continuous mapping.}
and hence
\begin{subequations}
\begin{equation}
\label{lambda_induction-T}
\va{\lambda}_T \preceq \va{x}_T.
\end{equation}
Then using \eqref{eqPMS1-b}, the optimality condition \eqref{eqHPMS1-eq1}
takes the form:
\begin{equation*}
\bespc{(H_{T-1})'_{u}
       \big(\va{x}_{T-1},\va{u}_{T-1},\va{w}_{T},
            \mu_{T-1} \compo f_{T-1}(\va{x}_{T-1},\va{u}_{T-1},\va{w}_{T})\big)}
      {\tribu{F}_{T-1}} = 0 \eqfinp
\end{equation*}
$\va{X}_{T-1}$ and $\va{u}_{T-1}$ being both $\tribu{F}_{T-1}$-measurable
random variables, and $\va{w}_{T}$ being independent of $\tribu{F}_{T-1}$
(white noise assumption), we deduce that the conditional expectation in
the last expression reduces to an expectation. Let $G_{T-1}$ denotes the
function resulting from its integration, namely
\begin{equation*}
G_{T-1}(x,u) = \Besp{(H_{T-1})'_{u}
\big(x,u,\va{w}_{T},\mu_{T-1} \compo f_{T-1}(x,u,\va{w}_{T})\big)}
\eqfinp
\end{equation*}
$G_{T-1}$ is a mesurable mapping,\footnote{in fact a continuous one
(from Assumption \ref{AssDiff})} and the optimality condition writes
\begin{equation*}
G_{T-1} (\va{x}_{T-1},\va{u}_{T-1}) = 0 \eqfinp
\end{equation*}
Using the measurable selection theorem available for implicit measurable
functions \cite[Theorem~8]{Leese:1974},\footnote{See for instance
\cite[Section~7]{Wagner:1977} for a survey of measurable selection theorems
corresponding to the implicit case. Note that \cite[Theorem~8]{Leese:1974}
needs a particular assumption concerning the $\sigma$-field equipping
$\espace{X}_{T-1}$. We assume here that such assumption holds.}
we deduce that there exists a measurable mapping
$\gamma_{T-1} : \espacea{X}_{T-1} \rightarrow \espacea{U}_{T-1}$
such that $G_{T-1} \big(\va{x}_{T-1},\gamma_{T-1}(\va{x}_{T-1})\big)=0$.
As a conclusion, the control variable
$\va{u}_{T-1}=\gamma_{T-1}(\va{x}_{T-1})$ satisfies the optimality
condition \eqref{eqHPMS1-eq1} at $t=T-1$ and is such that
\begin{equation}
\label{u_induction-T}
\va{u}_{T-1} \preceq \va{x}_{T-1}.
\end{equation}
\end{subequations}
\item At stage $t$, assume that
$\va{\lambda}_{t+1} \preceq (\va{x}_{t+1},\va{w}_{t+2},\ldots,\va{w}_T)$.
Then there exists a measurable function $\mu_{t}$ such that
\begin{subequations}
\begin{equation}
\label{lambda_induction-t}
\va{\lambda}_{t+1} = \mu_{t} (\va{x}_{t+1},\va{w}_{t+2},\ldots,\va{w}_T)
\eqfinp
\end{equation}
The optimality condition \eqref{eqHPMS1-eq1} at stage $t$ takes the form
\begin{equation*}
\Bespc{(H_{t})'_{u}
       \Big(\va{x}_{t},\va{u}_{t},\va{w}_{t+1},
            \mu_{t} \big(f_{t}(\va{x}_{t},\va{u}_{t},\va{w}_{t+1}),
            \va{w}_{t+2},\ldots,\va{w}_{T}\big)\Big)}
      {\tribu{F}_{t}} = 0 \eqfinp
\end{equation*}
With the same reasoning as at stage $T$, we deduce that this conditional
expectation reduces to an expectation, so that the optimality condition writes
\begin{equation*}
G_{t} (\va{x}_{t},\va{u}_{t}) = 0 \eqfinv
\end{equation*}
$G_{t}$ being a measurable function given by
\begin{equation*}
G_{t}(x,u) = \Besp{(H_{t})'_{u}
\Big(x,u,\va{w}_{t+1},\mu_{t} \big( f_{t}(x,u,\va{w}_{t+1}),
\va{w}_{t+2},\ldots,\va{w}_{T}\big)\Big)}
\eqfinp
\end{equation*}
Using again \cite[Theorem~8]{Leese:1974},
we deduce that there exists a measurable mapping
$\gamma_{t} : \espacea{X}_{t} \rightarrow \espacea{U}_{t}$
such that $\va{u}_{t}=\gamma_{t}(\va{x}_{t})$ satisfies
the optimality condition \eqref{eqHPMS1-eq1} at $t$.
We have accordingly
\begin{equation}
\label{u_induction-t}
\va{u}_{t} \preceq \va{x}_{t}.
\end{equation}
Ultimately, starting from the optimality condition \eqref{eqHPMS1-d},
namely
\begin{equation*}
\va{\lambda}_{t} = (H_{t})'^\top_{x}
(\va{x}_{t},\va{u}_{t},\va{w}_{t+1},\va{\lambda}_{t+1})
\eqfinv
\end{equation*}
using the induction assumption \eqref{lambda_induction-t}
together with \eqref{u_induction-t} and \eqref{eqPMS1-c}, we obtain that
\begin{equation*}
\va{\lambda}_{t} = (H_{t})'^\top_{x}
\Big(
  \va{x}_{t},\gamma_{t}(\va{x}_{t}),\va{w}_{t+1},
  \mu_{t}\big(f_{t}(\va{x}_{t},\gamma_{t}(\va{x}_{t}),\va{w}_{t+1}),
              \va{w}_{t+2},\ldots,\va{w}_T\big)
\Big)
\eqfinp
\end{equation*}
We conclude that
$\va{\lambda}_{t}\preceq\big(\va{x}_{t},\va{w}_{t+1},\ldots,\va{w}_T\big)$
so that the desired result holds true.\footnote{Note that we obtained as
an intermediate result that
$\va{\lambda}_{t+1}\preceq\big(\va{x}_{t},\va{w}_{t+1},\ldots,\va{w}_T\big)$.}
\end{subequations}
\end{itemize}

Let's go now to the general case $\Uas_t\subset\espacef{U}$.
From \eqref{cond123},  the optimality condition \eqref{eqHPMS1-e}
is again equivalent to an equality condition:
\begin{equation*}
\proj{\Uas_t}
     {\va{u}_{t}-\epsilon
      \bespc{{(H_{t})'_{u}}^{\top}
             (\va{x}_{t},\va{u}_{t},\va{w}_{t+1},\va{\lambda}_{t+1})}
            {\tribu{F}_{t}}} - \va{u}_{t} = 0 \eqfinv
\end{equation*}
and the same arguments as in the previous case remain valid.

\medskip

At last, from $\va{u}_t\preceq\va{x}_t$,
$\va{\lambda}_{t+1} \preceq (\va{x}_{t},\va{w}_{t+1},\ldots,\va{w}_T)$
and the white noise assumption,
we deduce that $\va{\lambda}_{t+1}$ depends on $\tribu{F}_t$ only through
$\va{x}_t$, so that
\begin{multline*}
\bespc{(L_t)'_u (\va{x}_t,\va{u}_t,\va{w}_{t+1})+
         \va{\lambda}_{t+1}^\top (f_t)'_u (\va{x}_t,\va{u}_t,\va{w}_{t+1})}
       {\tribu{F}_t}= \\
\bespc{(L_t)'_u (\va{x}_t,\va{u}_t,\va{w}_{t+1})+
         \va{\lambda}_{t+1}^\top (f_t)'_u (\va{x}_t,\va{u}_t,\va{w}_{t+1})}
       {\va{x}_t} \eqfinp
\end{multline*}
Statement~$(c)$ is thus satisfied and the proof is complete. \qquad
\end{proof}

In the Markovian case, the optimal solution of Problem~\eqref{SOC}
(measurability with respect to all past noise variables)
satisfies the measurability constraints of Problem~\eqref{SOCX}
(measurability with respect to the current state variable).
The two problems are equivalent (same $\min$ and $\argmin$).
In fact, the feasible set of~\eqref{SOC} contains the feasible set
of~\eqref{SOCX}, and Theorem~\ref{thmPMS1X} shows us that any
optimal solution of~\eqref{SOC} is feasible for~\eqref{SOCX},
and is therefore also optimal also for~\eqref{SOCX}.
Ultimately, the optimality conditions of Problem~\eqref{SOCX}
can be written as
\begin{subequations}
\begin{align}
& \va{x}_0=\va{w}_0 \eqfinv \\
& \va{x}_{t+1} = f_{t}(\va{x}_{t},\va{u}_{t},\va{w}_{t+1}) \eqfinv \\
& \nonumber \\
& \va{\lambda}_{T} = K'^\top(\va{x}_{T}) \eqfinv \\
& \va{\lambda}_{t} = (L_{t})'^\top_{x_t}(\va{x}_{t},\va{u}_{t},\va{w}_{t+1})+
                     (f_{t})'^\top_{x_t}(\va{x}_{t},\va{u}_{t},\va{w}_{t+1})
                     \va{\lambda}_{t+1} \eqfinv \\
& \nonumber \\
& \bespc{(L_{t})'_{u}(\va{x}_{t},\va{u}_{t},\va{w}_{t\!+\!1})\!+\!
          \va{\lambda}_{t\!+\!1}^\top(f_{t})'_{u}(\va{x}_{t},\va{u}_{t},\va{w}_{t\!+\!1})}
        {\va{x}_{t}} \in -\partial\fcara{\Uas_t}(\va{u}_t) \eqfinp
\end{align}
\label{eqPMS1X}
\end{subequations}

\begin{remark}
Let $(\va{g}_t^1)_{t=0,\ldots,T-1}$ and $(\va{g}_t^2)_{t=0,\ldots,T-1}$ be
the gradient  processes associated with~\eqref{eqPMS1} and~\eqref{eqPMS1X}
respectively:
\begin{align*}
\va{g}_t^1\defegal &
\bespc{(L_{t})'^\top_{u}(\va{x}_{t},\va{u}_{t},\va{w}_{t+1})+
       (f_{t})'^\top_{u}(\va{x}_{t},\va{u}_{t},\va{w}_{t+1})\va{\lambda}_{t+1}}
      {\tribu{F}_{t}} \eqfinv \\
\va{g}_t^2\defegal &
\bespc{(L_{t})'^\top_{u}(\va{x}_{t},\va{u}_{t},\va{w}_{t+1})+
       (f_{t})'^\top_{u}(\va{x}_{t},\va{u}_{t},\va{w}_{t+1})\va{\lambda}_{t+1}}
      {\va{x}_{t}} \eqfinp
\end{align*}
Unlike Problem \ref{SOC}, we are unable to compute the optimality
conditions~\eqref{eqPMS1X} of Problem~\eqref{SOCX} by differentiating
the Lagrangian function, because the conditioning term $\va{x}_{t}$
depends itself on the control variables $\va{u}_{t}$.
Consequently, $\va{g}_{t}^{2}$ is not claimed to represent the projected
gradient of Problem~\eqref{SOCX}. The equality between the
gradient $\va{g}_t^1$ and $\va{g}_t^2$ holds true only at the optimum.
\end{remark}

\subsubsection{Markovian case: adapted optimality conditions}

We now present the adapted version of the optimality
conditions of Problem~\eqref{SOC} under Markovian assumptions.

%\vspace{0.2cm}

%\begin{assumption}[Measurable selection]
%\begin{enumerate}
%\item[] For any $t=0,\ldots,T-1$, any $x\in\espacea{X}_t$ and any
%$(\va{\Lambda}_t)_{t=0,\ldots,T}$ solution of~\eqref{eqPMS2},
%the mapping $g_t(x,u)\defegal\besp{(L_{t})'_u(x,u,\va{w}_{t+1})+
%             \va{\Lambda}_{t+1}^\top(f_{t})'_u(x,u,\va{w}_{t+1})}$
%is continuous with respect to $u$, and $\exists \: u\in\espacea{U}_t$
%such that $g_t(x,u)\in-\partial\fcara{\Gammaas_t}(u)$.
%\end{enumerate}
%\label{AssMeaSel2}
%\end{assumption}

\vspace{0.2cm}

\begin{theorem}
Suppose that Assumptions~\ref{AssStruct}, \ref{AssDiff} and \ref{AssMark}
are fulfilled, and assume that there exists two
random processes $(\va{u}_t)_{t=0,\ldots,T-1}\in\espacef{U}$ and
$(\va{x}_t)_{t=0,\ldots,T-1}\in\espacef{X}$ solution of Problem
\eqref{SOC}. Then there exists a process
$(\va{\Lambda}_t)_{t=0,\ldots,T-1}\in\espacef{X}$ satisfying
\eqref{eqPMS2} and such that, for all $t=0,\ldots,T-1$,
\begin{align*}
& (a) \enspace
  \va{\Lambda}_{t+1}\preceq \va{x}_{t+1} \eqfinv \\
& (b) \enspace
  \va{u}_t \preceq \va{x}_t \eqfinv \\
& (c) \enspace
  \bespc{(L_t)'_u(\va{x}_t,\va{u}_t,\va{w}_{t+1})+
          \va{\Lambda}_{t+1}^\top(f_t)'_u(\va{x}_t,\va{u}_t,\va{w}_{t+1})}
        {\tribu{F}_t} = \\
& \quad\quad\quad\quad\quad
  \bespc{(L_t)'_u(\va{x}_t,\va{u}_t,\va{w}_{t+1})+
         \va{\Lambda}_{t+1}^\top(f_t)'_u(\va{x}_t,\va{u}_t,\va{w}_{t+1})}
        {\va{x}_t} \eqsepv \Pps \eqfinp
\end{align*}
\label{thmPMS2X}
\end{theorem}

\begin{proof}
The proof follows the same scheme as the one of Theorem~\ref{thmPMS1X}.
We just point out that
$\va{\lambda}_{t+1}\preceq (\va{x}_{t+1},\va{w}_{t+2},\ldots,\va{w}_T)$
and $\va{\Lambda}_{t+1}=\nespc{\va{\lambda}_{t+1}}{\tribu{F}_{t+1}}$
implies that $\va{\Lambda}_{t+1}\preceq\va{x}_{t+1}$. \qquad
\end{proof}

From the equivalence between Problem~\eqref{SOC} (measurability with
respect to the past noise) and Problem~\eqref{SOCX} (measurability with
respect to the state) in the Markovian framework, we may consider that
the optimality conditions of Problem~\eqref{SOCX}  in the adapted version
are
\begin{subequations}
\begin{align}
& \va{x}_0=\va{w}_0 \eqfinv \\
& \va{x}_{t+1} = f_{t}(\va{x}_{t},\va{u}_{t},\va{w}_{t+1}) \eqfinv \\
& \nonumber \\
& \va{\Lambda}_{T} = K'^\top(\va{x}_{T}) \eqfinv \\
& \va{\Lambda}_{t} = \bespc{(L_{t})'^\top_{x}(\va{x}_{t},\va{u}_{t},\va{w}_{t+1})+
                            (f_{t})'^\top_{x}(\va{x}_{t},\va{u}_{t},\va{w}_{t+1})
                            \va{\Lambda}_{t+1}}
                           {\va{x}_t} \eqfinv \\
& \nonumber \\
& \bespc{(L_{t})'_{u}(\va{x}_{t},\va{u}_{t},\va{w}_{t\!+\!1})\!+\!
          \va{\Lambda}_{t\!+\!1}^\top(f_{t})'_{u}(\va{x}_{t},\va{u}_{t},\va{w}_{t\!+\!1})}
        {\va{x}_{t}} \in -\partial\fcara{\Uas_t}(\va{u}_t) \eqfinp
\end{align}
\label{eqPMS2X}
\end{subequations}

\medskip

The optimality conditions~\eqref{eqPMS1X} and~\eqref{eqPMS2X} involve
conditional expectations with respect to the state variable
the dimension of which is, in most cases, fixed, that is, it does
not depend on the time stage. In order to solve Problem~\eqref{SOC}
(and equivalently Problem~\eqref{SOCX}), we have to discretize
those conditions and in particular to approximate the conditional
expectations. The literature on conditional expectation approximation
only offers biased estimators with an integrated squared
error depending on the dimension of the conditioning term.
On the contrary, approximating an expectation through a Monte-Carlo
technique involves non-biased estimators the variance of which does not
depend on the dimension of the underlining space. In the next section,
we propose a functional interpretation of the stochastic optimal
control problem in order to get rid of those conditional expectations
and deal only with expectations.

\subsection{\label{OCfunc}Optimality conditions from a functional point of view}

Consider the stochastic optimal control Problem~\eqref{SOC}.
Under Markovian assumptions, we have shown in \S\ref{KKTX}
that it is equivalent to Problem~\eqref{SOCX} and that~\eqref{eqPMS2X}
is a set of necessary optimality conditions.
Hereafter we transform the optimality conditions~\eqref{eqPMS2X}
using Theorem~\ref{thmPMS2X} and the functional interpretation
of the measurability relation between random variables
(see Proposition~\ref{propmes}).

\begin{theorem}
Suppose that Assumptions~\ref{AssStruct}, \ref{AssDiff} and \ref{AssMark}
are fulfilled, and assume that there exist two
random processes $(\va{u}_t)_{t=0,\ldots,T-1}\in\espacef{U}$ and
$(\va{x}_t)_{t=0,\ldots,T-1}\in\espacef{X}$ solution of Problem
\eqref{SOC}. Let $(\va{\Lambda}_t)_{t=0,\ldots,T-1}\in\espacef{X}$
be a random process satisfying the optimality conditions~\eqref{eqPMS2}.
Then there exists two sequences of mappings
$(\Lambda_t)_{t=0,\ldots,T}$ and $(\phi_t)_{t=0,\ldots,T-1}$,
$\Lambda_t:\espacea{X}_t\rightarrow\espacea{X}_t$ and
$\phi_t:\espacea{X}_t\rightarrow\espacea{U}_t$,
such that for all $t=0,\ldots,T-1$,
{\small
\begin{subequations}
\label{eqPMS3X}
\begin{align}
& \va{\Lambda}_t=\Lambda_t(\va{x}_t) \eqfinv \\
& \va{u}_t=\phi_t(\va{x}_t) \eqfinv \\
& \nonumber \\
& \Lambda_{t}(\cdot) =
  \Besp{(L_{t})'^\top_{x}\big(\cdot,\phi_{t}(\cdot),\va{w}_{t+1}\big)+
        (f_{t})'^\top_{x}\big(\cdot,\phi_{t}(\cdot),\va{w}_{t+1}\big)
        \Lambda_{t+1}\big(f_t(\cdot,\phi_t(\cdot),\va{w}_{t+1})\big)} \eqfinv \\
& \Besp{(L_{t})'_{u}\big(\cdot,\phi_{t}(\cdot),\va{w}_{t\!+\!1}\big)+
        \Lambda_{t+1}^\top\big(f_t(\cdot,\phi_t(\cdot),\va{w}_{t\!+\!1})\big)
        (f_{t})'_{u}\big(\cdot,\phi_{t}(\cdot),\va{w}_{t\!+\!1}\big)}
\!\in\!-\partial\fcara{\Phias_t}(\phi_t) \eqfinv
\end{align}
\end{subequations}
}

\noindent
with
$\Phias_t\defegal\left\{\phi_t\in
  L^2(\espacea{X}_t,\borel{\espacea{X}_t},\prbt_{\va{x}_t};\espacea{U}_t)
  \eqsepv \phi_t(x)\in\Gammaas_t \eqsepv \forall x\in\espacea{X}_t\right\}$
and $\prbt_{\va{x}_t}$ the probability measure associated with $\va{x}_t$.

\label{thmPMS3X}
\end{theorem}

\begin{proof}
By Theorem~\ref{thmPMS2X}, the random processes $(\va{u}_t)_{t=0,\ldots,T-1}$
$(\va{x}_t)_{t=0,\ldots,T}$ and $(\va{\Lambda}_t)_{t=0,\ldots,T-1}$ are such
that $\va{u}_t\preceq\va{x}_t$ and $\va{\Lambda}_t\preceq\va{x}_t$.
By Proposition~\ref{propmes}, there exist measurable mappings
$\phi_t:\espacef{X}_t\rightarrow\espacea{U}_t$ and
$\Lambda_t:\espacea{X}_t\rightarrow\espacea{X}_t$ such that
$\va{u}_t=\phi_t(\va{x}_t)$ and $\va{\Lambda}_t=\Lambda_t(\va{x}_t)$.
From $\va{u}_t\in\espacef{U}_t=L^2(\omeg,\trib,\prbt;\espacea{U}_t)$,
we deduce that
$\phi_t\in L^2(\espacea{X}_t,\borel{\espacea{X}_t},\prbt_{\va{x}_t};\espacea{U}_t)$:
\begin{equation*}
\int_\Omega\norm{\va{u}_t(\omega)}^2\mathrm{d}\prbt(\omega) =
\int_\Omega\norm{\phi_t\left(\va{x}_t(\omega)\right)}^2\mathrm{d}\prbt(\omega) =
\int_{\espacea{X}_t}\norm{\phi_t(x)}^2\mathrm{d}\prbt_{\va{x}_t}(x) <
+\infty.
\end{equation*}
Since $\va{u}_t\in\Uas_t \Leftrightarrow
\phi_t\left(\va{x}_t(\omega)\right)\in\Gammaas_t \eqsepv \Pps$,
we obtain that $\phi_t\in\Phias_t$.

The co-state dynamic equations in~\eqref{eqPMS2X} rewrites
{\small
\begin{equation*}
\Lambda_{t}(\va{x}_t) =
\Bespc{(L_{t})'^\top_{x}\big(\va{x}_{t},\phi_{t}(\va{x}_t),\va{w}_{t+1}\big)+ \\
       (f_{t})'^\top_{x}\big(\va{x}_{t},\phi_{t}(\va{x}_t),\va{w}_{t+1}\big)
       \Lambda_{t+1}(\va{x}_{t+1})}
      {\va{x}_t} \eqfinv
\end{equation*}
}
which, using $\va{x}_{t+1}=f_t\left(\va{x}_t,\phi_t(\va{x}_t),\va{w}_{t+1}\right)$,
becomes
{\small
\begin{multline*}
\Lambda_{t}(\va{x}_t) =
\Bespc{(L_{t})'^\top_{x}\big(\va{x}_t,\phi_{t}(\va{x}_t),\va{w}_{t+1}\big)\\
       +(f_{t})'^\top_{x}\big(\va{x}_t,\phi_{t}(\va{x}_t),\va{w}_{t+1}\big)
       \Lambda_{t+1}\big(f_t(\va{x}_t,\phi_t(\va{x}_t),\va{w}_{t+1})\big)}
      {\va{x}_t} \eqfinp
\end{multline*}
}
From the white noise assumption, the random variable $\va{w}_{t+1}$ is
independent of $\va{x}_t$, so that the conditional expectation turns
out to be just an expectation with respect to the noise random variable.
The co-state dynamics equation writes accordingly as a functional equality:
{\small
\begin{equation*}
%%%%% \Lambda_{T}(\cdot) = K'^\top(\cdot),\\
\Lambda_{t}(\cdot) =
\Besp{(L_{t})'^\top_{x}(\cdot,\phi_{t}(\cdot),\va{w}_{t+1})+
      (f_{t})'^\top_{x}(\cdot,\phi_{t}(\cdot),\va{w}_{t+1})
      \Lambda_{t+1}\big(f_t(\cdot,\phi_t(\cdot),\va{w}_{t+1})\big)}.
\end{equation*}
}

Using similar arguments, we easily obtain from the last condition
in~\eqref{eqPMS2X}:
{\small
\begin{equation*}
\Besp{(L_{t})'_{u}(\cdot,\phi_{t}(\cdot),\va{w}_{t+1})+
      \Lambda_{t+1}^\top\big(f_t(\cdot,\phi_t(\cdot),\va{w}_{t+1})\big)
      (f_{t})'_{u}(\cdot,\phi_{t}(\cdot),\va{w}_{t+1})}
\in-\partial\fcara{\Phias_t}(\phi_t) \eqfinv
\end{equation*}
}
which completes the proof. \qquad
\end{proof}

Theorem~\ref{thmPMS3X} provides the new functional optimality
conditions~\eqref{eqPMS3X} for Problem~\eqref{SOC} in the Markovian
case. These optimality conditions do not involve conditional
expectations but just expectations. Therefore, we may hope
that, in the approximation process of these conditions, we will
obtain non biased estimates the variance of which will not depend on
the dimension of the state space.

\section{\label{adapdisc}Adaptive discretization technique}

In this section we develop \emph{tractable} numerical methods
for obtaining the solution of Problem~\eqref{SOC}. We will limit
ourselves here to the Markovian framework, but methods for all
cases described in \S\ref{OCalg} can be found in \cite{TheseDallagi}.
%Note however that, at least from the theoretical point of view, the
%Markovian assumption could be considered as not restrictive. As a matter
%of fact, if the noises are the output of an ARMA process, one should increase
%the vector $\va{x}_t$ to account for the memory of the noise process.
%Nevertheless such a modelling introduces computational difficulties:
%increasing the state dimension often means giving up using dynamic
%programming. This will be explained later in \S\ref{adapdisc:dynprog}.
We briefly discuss two classical solution methods, namely
\emph{stochastic programming} and \emph{dynamic programming}.
Then we present an \emph{adaptive mesh} algorithm which consists in
discretizing the optimality conditions obtained at \S\ref{OCfunc}
and in using them in a gradient-like algorithm.

\subsection{Discrete representation of a function\label{RmqRT}}
As far as numerical resolution is concerned, we need to manipulate
functions which are infinite dimensional objects and which, in most
cases, do not have a closed-form expression. Thus we must have
a discrete representation of such an object.
Let $\phi : \espacea{X}\rightarrow\espacea{U}$ be a function.
We suppose that we have at disposal a fixed or variable grid
$\grid{x}$ in $\espacea{X}$, that is a collection of elements
in $\espacea{X}$:
\begin{equation*}
\grid{x} = (x^i)_{i=1,\ldots,n} \in \espacea{X}^n \eqfinp
\end{equation*}
A point $x^{i}$ will be called a \emph{particle}. In order to obtain a discrete
representation of $\phi$, we define its \emph{trace}, that is a grid $\grid{u}$
on $\espacea{U}$ which corresponds to the values of $\phi$ on $\grid{x}$:
\begin{equation*}
\grid{u} = \big(\phi(x^i)\big)_{i=1,\ldots,n} \in \espacea{U}^n \eqfinp
\end{equation*}
We denote by $\gridop{T}_{\espacea{U}}:
\espacea{U}^\espacea{X}\rightarrow\espacea{X}^n\times\espacea{U}^n$
the \emph{trace operator} which, with any function
$\phi:\espacea{X}\rightarrow\espacea{U}$, associates the couple
of grids $(\grid{x},\grid{u})$, that is the $n$ points
$\big(x^i,\phi(x^i)\big)\in\espacea{X}\times\espacea{U}$.
On the other hand, knowing the trace of $\phi$, we need to compute
an approximation of the value taken by $\phi$ at any $x\in\espacea{X}$.
Let $\gridop{R}_{\espacea{U}}:
 \espacea{X}^n\times\espacea{U}^n\rightarrow\espacea{U}^\espacea{X}$
be an \emph{interpolation-regression operator} which, with any couple
of grids $(\grid{x},\grid{u})$, associates
$\widetilde{\phi}:\espacea{X}\rightarrow\espacea{U}$ representing the initial
function. Such an interpolation-regression operator may be defined in different
ways (polynomial interpolation, kernel approximation, closest neighbor, etc).

\subsection{Stochastic programming\label{adapdisc:stochprog}}

One way to discretize stochastic optimal control problems of type
\eqref{SOC} is to model the information structure as a decision
tree.
\begin{enumerate}
\item First, simulate a given number $N$ of scenarios
$(\va{w}_t^k)_{t=0,\ldots,T}^{k=1,\ldots,N}$ of the noise process.
Then, by some tree generation procedures (see e.g. \cite{Pflug01} or
\cite{HeitschRomisch03}), organize these scenarios in a scenario tree
(at any node of any time stage $t$, there is a single past trajectory
but multiple futures: see Figure~\ref{treegen}).
\begin{figure}[htbp]
\begin{center}
\includegraphics[width=10cm,clip]{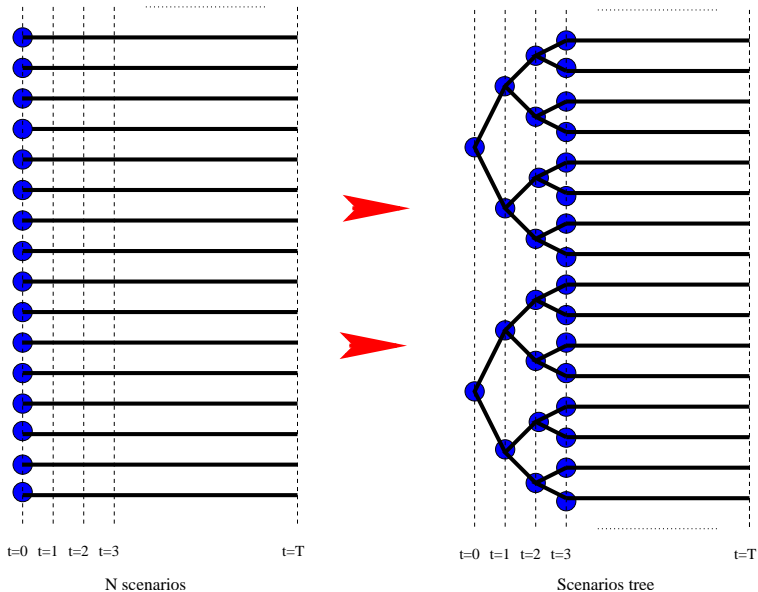}
\caption{\label{treegen}From scenarios to a tree}
\end{center}
\end{figure}
\item Second, write the components of the problem (state dynamics
and cost expectation) on  the scenario tree (note that the information
constraints are built-in in such a tree structure).
\end{enumerate}
Then the approximation on the scenario tree of Problem~\eqref{SOC}
is solved using an appropriate (deterministic) non-linear
programming package. The optimal solution consists of state and control
particles at each node of the tree. An interpolation-regression procedure
(as suggested at~\S\ref{RmqRT}) has to be performed at each time stage
in order to synthesize a feedback law.

Such a methodology is relatively easy to implement and need not in fact
any Markovian assumption. Nevertheless, it faces a serious drawback:
at the first time stages (close to the tree root), only a few particles
--- or nodes --- are available, and experience shows that the estimates
of the optimal feedback they provide has a relatively weak variance
(because the number of pending scenarios at each such node is still
large enough); on the contrary, at the final time stages (near the tree
leaves), a large number of particles is available, but with a huge variance.
In both cases, the feedback synthesis (the interpolation-regression process
from the available grid) will be inaccurate. Such an observation will be
highlighted later on the case study~\S\ref{SecTest}.

\subsection{Stochastic dynamic programming\label{adapdisc:dynprog}}

We consider Problem~\eqref{SOC} in the Markovian case, which
is thus equivalent to Problem~\eqref{SOCX}. From Theorem~\ref{thmPMS1X},
the optimal control process $(\va{u}_t)_{t=0,\ldots,T-1}$
can be searched as a collection of feedback laws
$(\phi_t)_{t=0,\ldots,T-1}$ depending on the process
state $(\va{x}_t)_{t=0,\ldots,T}$, that is
\begin{equation*}
\va{u}_t=\phi_t(\va{x}_t)
\eqsepv \Pps \eqfinp
\end{equation*}

According to the \emph{Dynamic Programming Principle}, the resolution
is built up \emph{backward} from $t=T$ to $t=0$ by solving the Bellman equation at each time
stage for all state values $x\in\espacea{X}_t$.
Such a principle leads to the following algorithm
\cite{Bertsekas:1976,Bertsekas-Shreve:1996}.

\begin{algorithm}[Stochastic Dynamic Programming]
\begin{itemize}
\item \textbf{At stage $T$},\\
define the Bellman function $V_{T}$ as
\begin{equation*}
V_T(x)\defegal K(x) \eqsepv \forall  x\in\espacea{X}_T \eqfinp
\end{equation*}
\item \textbf{Recursively for $t=T-1,\ldots,0$},\\
compute the Bellman function $V_{t}$ as
\begin{equation*}
V_t(x) = \min_{u\in \Gammaas_t}
\Besp{L_t(x,u,\va{w}_{t+1})+V_{t+1}\big(f_t(x,u,\va{w}_{t+1})\big)}
\eqsepv \forall x\in\espacea{X}_t \eqfinv
\end{equation*}
the optimal feedback law $\phi_{t}$ being obtained as
\begin{equation*}
\phi_t(x)=\mathrm{arg}\min_{u\in \Gammaas_t}
\Besp{L_t(x,u,\va{w}_{t+1})+V_{t+1}\big(f_t(x,u,\va{w}_{t+1})\big)}
\eqsepv \forall x\in\espacea{X}_t \eqfinp
\end{equation*}
\end{itemize}
\label{algoDP}
\end{algorithm}

This algorithm is only conceptual because it operates on infinite
dimensional objects $V_t$ (and expectations cannot always be evaluated
analytically). We must indeed manipulate those objects as indicated
at~\S\ref{RmqRT}. For every $t=0,\ldots,T$, let
$\grid{x}_t\defegal (x_t^i)_{i=1,\ldots,n_t}$ be a \emph{fixed} grid
of $n_t$ discretization points in the state space $\espacea{X}_t$.
We approximate the functions appearing in Algorithm~\ref{algoDP}
by their trace over that grid:
\begin{alignat*}{3}
& v_t^i= V_t(x_t^i) \eqsepv
   & \forall i=1,\ldots,n_t \eqsepv
   & \forall t=0,\ldots,T \eqfinv \\
& u_t^i= \phi_t(x_t^i) \eqsepv
   & \forall i=1,\ldots,n_t \eqsepv
   & \forall t=0,\ldots,T-1 \eqfinp
\end{alignat*}
We also need to approximate the expectations by the Monte Carlo method.
Let $(\va{w}_t^k)_{t=0,\ldots,T}^{k=1,\ldots,N}$ denote $N$ independent
and identically distributed scenarios of the random noise process.
The discretized stochastic dynamic programming algorithm is as follows.

\begin{algorithm}[Discretized Stochastic Dynamic Programming]
\begin{itemize}
\item \textbf{At stage $T$},\\
compute the trace $\grid{v}_{T}$ of the Bellman function $V_{T}$:
\begin{equation*}
v_T^i=V_T(x_T^i)\eqsepv \forall i=1,\ldots,n_T \eqfinv
\end{equation*}
\item \textbf{Recursively for $t=T-1,\ldots,0$},\\
approximate $V_{t+1}$ by interpolation-regression:
\begin{equation*}
\widehat{V}_{t+1}=\gridop{R}_\bbR\left(\grid{x}_{t+1},\grid{v}_{t+1}\right) \eqfinv
\end{equation*}
compute the two grids $\grid{v}_{t}$ and $\grid{u}_{t}$, that is,
for each $i=1,\ldots,n_t$,
\begin{align*}
v_t^i & =
   \min_{u\in \Gammaas_t}
   \frac{1}{N}\sum_{k=1}^N\biggr[ L_t\big(x_t^i,u,\va{w}_{t+1}^k\big)+
   \widehat{V}_{t+1}\big(f_t(x_t^i,u,\va{w}_{t+1}^k)\big) \biggr] \eqfinv \\
u_t^i & =
   \argmin_{u\in \Gammaas_t}
   \frac{1}{N}\sum_{k=1}^N\biggr[ L_t\big(x_t^i,u,\va{w}_{t+1}^k\big)+
             \widehat{V}_{t+1}\big(f_t(x_t^i,u,\va{w}_{t+1}^k)\big) \biggr] \eqfinv
\end{align*}
and obtain the feedback law as
\begin{equation*}
\phi_t=\gridop{R}_{\espacea{U}_t}\left(\grid{x}_t,\grid{u}_t\right) .
\end{equation*}
\end{itemize}
\label{DPApprox}
\end{algorithm}

\begin{remark}
The interpolation-regression operator is in most cases mandatory
in Algorithm \ref{DPApprox}. As a matter of fact, for a given time stage
$t$ and a given control value $u\in\espacea{U}_t$, there usually does not
exist any index $j$ such that $x_{t+1}^j= f_t(x_t^i,u,\va{w}_{t+1}^k)$,
which means that $V_{t+1}$ has to be computed out of the grid $\grid{x}_{t+1}$.
Note however that the Bellman function at time stage $T$ is known analytically,
so that interpolation is not needed for $V_{T}$.
\end{remark}

This method faces an important difficulty: the \emph{curse of dimensionality}.
In fact, one generally discretizes each state coordinate at each time
stage $t$ using a \emph{scalar} grid and a fixed number of points. Therefore,
the grid $\grid{x}_t$ is obtained as the Cartesian product of the scalar grids
over all the coordinates. Thus, the number of particles of that grid increases
\emph{exponentially} with the state space dimension. This is the well-known
drawback of most methods derived from Dynamic Programming, which do not take
advantage of the repartition of the optimal state particles in the state space
in order to concentrate computations in significant parts of the state space.

\subsection{The adaptive mesh algorithm}

Considering the difficulties faced by both stochastic and dynamic
programming methods, we propose an alternative method for solving
Problem~\eqref{SOC} in the Markovian case. The method, based on
optimality conditions~\eqref{eqPMS3X}, aims at
\begin{itemize}
\item dealing with the same number of noise particles from the beginning
of the time horizon to the end: we thus hope that the generated feedback
law estimators will have a reduced and fixed variance during all time stages;
\item attempting to alleviate the curse of dimensionality by operating
on an adaptive discretization grid automatically generated from the primary
noise discretization grid.
\end{itemize}

\subsubsection{Approximation}

Let us denote by $(\va{w}_t^k)_{t=0,\ldots,T}^{k=1,\ldots,N}$ a set
of $N$ independent and identically distributed scenarios obtained
from the noise random process. Given random control grids
$\grid{u}_t\defegal (\va{u}_t^k)_{k=1,\ldots,N}$ for $t=0,\ldots,T-1$,
we can compute the state random grids
$\grid{x}_t\defegal (\va{x}_t^k)_{k=1,\ldots,N}$
by propagating the state dynamics equation:
\begin{subequations}
\begin{align}
& \va{x}_0^k=\va{w}_0^k \eqsepv \forall k=1,\ldots,N \eqfinv \\
& \va{x}_{t+1}^k=f_t\left(\va{x}_t^k,\va{u}_t^k,\va{w}_{t+1}^k\right)
  \eqsepv \forall k=1,\ldots,N \eqsepv \forall t=0,\ldots,T-1 \eqfinp
\end{align}
\label{DynXApprox}
\end{subequations}
The feedback at time stage $t$ is obtained using an interpolation-regression
operator on the grids $(\grid{x}_t,\grid{u}_t)$. Note that the state space
grids $\grid{x}_t$ are not a priori fixed as in Dynamic Programming methods.
They are in fact adapted to the control particules, as they change whenever
the decision maker changes its control strategy.

\begin{remark}
If the optimal state is concentrated in a region of the state space,
the optimal feedback will be synthesized only inside that region.
In fact we need not compute it elsewhere, because the state will hardly
reach other regions.
\end{remark}

In order to obtain approximations $\widehat{\Lambda}_t$ of the co-state
functions $\Lambda_t$ introduced at Theorem \ref{thmPMS3X}, we compute
particles $\va{\Lambda}^k_t$ by integrating the co-state backward dynamic
equation over the state grids $\grid{x}_t$. We thus obtain co-state grids
$\grid{l}_t=(\va{\Lambda}^k_t)_{k=1,\ldots,N}$, and make use of
interpolation-regression operators $\gridop{R}_{\espacea{X}_t}$
in order to compute the co-state function for values out of the current
grid. More specifically, the process is initiated with
$\widehat{\Lambda}_T(\cdot) = \Lambda_T(\cdot) = K'(\cdot)$;
then, for all $t=T-1,\ldots,0$, one computes:
{\small
\begin{subequations}
\begin{align}
& \va{\Lambda}^k_t = \frac{1}{N}\sum_{j=1}^N
     (L_t)'^\top_x(\va{x}_t^k,\va{u}_t^k,\va{w}_{t+1}^j)+
     (f_t)'^\top_x(\va{x}_t^k,\va{u}_t^k,\va{w}_{t+1}^j)
     \widehat{\Lambda}_{t+1}^\top\big(f_t(\va{x}_t^k,\va{u}_t^k,\va{w}_{t+1}^j)\big) \eqfinv \\
& \widehat{\Lambda}_t(\cdot)=\gridop{R}_{\espacea{X}_t}(\grid{x}_t,\grid{l}_t)(\cdot) \eqfinp
\end{align}
\label{DynLApprox}
\end{subequations}
}

Ultimately, for all $k=1,\ldots,N$ and for all $t=0,\ldots,T-1$, the gradient
particles $\va{g}_t^k$ are obtained as
{\small
\begin{equation}
\va{g}_t^k=\frac{1}{N}\sum_{j=1}^N\biggr[
(L_t)'^\top_u\big(\va{x}_t^k,\va{u}_t^k,\va{w}_{t+1}^j\big)+
(f_t)'^\top_u\big(\va{x}_t^k,\va{u}_t^k,\va{w}_{t+1}^j\big)
\widehat{\Lambda}_{t+1}^\top\Big(f_t\big(\va{x}_t^k,\va{u}_t^k,\va{w}_{t+1}^j\big)\Big)
\biggr] \eqfinp
\label{Gradapprox}
\end{equation}
}
As already noticed, the direction associated with these particles
represents the gradient only at the optimum.

\subsubsection{Algorithm}

We can now derive a descent-like algorithm to solve Problem \eqref{SOC}
under Markovian assumptions. At each iteration, state particles are
propagated forward --- with no interaction between particles --- then,
co-state particles are propagated backward --- now with interaction
caused by the regression-interpolation operations). Then, gradient
particles are computed using~\eqref{Gradapprox} and the control particles
are updated using a gradient-like  method. Ultimately, a functional
representation of the feedback laws is obtained thanks to
a regression-interpolation operator.

\begin{algorithm}~
\begin{itemize}
\item \textbf{Step $[0]$}.\\
Let $(\grid{u}_t^{[0]})_{t=0,\ldots,T-1}=
     \big(\va{u}_t^{[0],k}\big)_{t=0,\ldots,T-1}^{k=1,\ldots,N}$
be the initial control grids.
\item \textbf{Step $[\ell]$}.
\begin{enumerate}
\item Compute the state grids $(\grid{x}_t^{[\ell]})_{t=0,\ldots,T}$ by
propagating the dynamics~\eqref{DynXApprox} with $\va{u}=\va{u}^{[\ell]}$.
\item Compute both co-state grids $(\grid{l}_t^{[\ell]})_{t=0,\ldots,T}$ and
functional approximations $(\widehat{\Lambda}_{t}^{[\ell]})_{t=0,\ldots,T-1}$
by propagating the dynamics~\eqref{DynLApprox} with $\va{u}=\va{u}^{[\ell]}$
and $\va{x}=\va{x}^{[\ell]}$.
\item Compute the gradient particles
$\big(\va{G}_t^{[\ell],k}\big)_{t=0,\ldots,T-1}^{k=1,\ldots,N}$
using Equation~\eqref{Gradapprox}
with $\va{u}=\va{u}^{[\ell]}$, $\va{x}=\va{x}^{[\ell]}$ and
$\widehat{\Lambda}=\widehat{\Lambda}^{[\ell]}$.
\item For all $t=0,\ldots,T-1$ and $k=1,\ldots,N$, update the control particles
by performing a projected gradient step:
\begin{equation*}
\va{u}_t^{[\ell+1],k} =
\proj{\Gammaas_{t}}{\va{u}_t^{[\ell],k}-\rho^{[\ell]}\va{g}_{t}^{[\ell],k}} \eqfinp
\end{equation*}
\item Stop if some degree of accuracy is reached. Else set $\ell=\ell+1$ and iterate.
\end{enumerate}
\item \textbf{Step $[\infty]$}.
\begin{enumerate}
\item Set the grids:
\begin{align*}
(\grid{u}_t^{[\infty]})_{t=0,\ldots,T-1}&=
  \big(\va{u}_t^{[\ell],k}\big)_{t=0,\ldots,T-1}^{k=1,\ldots,N} \eqfinv \\
(\grid{x}_t^{[\infty]})_{t=0,\ldots,T-1}&=
  \big(\va{x}_t^{[\ell],k}\big)_{t=0,\ldots,T-1}^{k=1,\ldots,N} \eqfinp
\end{align*}
\item Obtain for all $t=0,\ldots,T-1$ the feedback law:
\begin{equation*}
\phi_t(\cdot)=
\gridop{R}_{\espacea{U}_t}
\big(\grid{x}^{[\infty]}_t,\grid{u}_t^{[\infty]}\big)(\cdot)
\eqfinp
\end{equation*}
\end{enumerate}
\end{itemize}
\label{gradproj}
\end{algorithm}

The only a priori discretization made during this algorithm is relative to
the noise sampling. Once such a discretization has been performed, all other grids
used to ultimately obtain the feedback laws are derived by integrating
dynamic equations. In addition, no conditional expectation approximations
are involved in the process. We just approximate expectations using
Monte Carlo techniques, and it is well-known that the variance of such
an approximation does not depend on the dimension of the underlying space.
The only space-size dependent operations are in fact the interpolation operators
used to approximate the co-state mappings during the iterative process, plus
the feedback laws once convergence is achieved.

Furthermore, this adaptive discretization makes computations concentrate
in effective state space regions, unlike dynamic programming which explores
the whole state space.

%Therefore, we hope that the approximation error generated by the adaptive
%mesh algorithm does not depend on the space dimension, thus avoiding
%the curse of dimensionality.

\section{\label{SecTest}Case study}

We consider the production management of an hydro-electric dam.
The problem is formulated as a stochastic optimal control,
and we consider solving it by the three methods described
in \S\ref{adapdisc}.

\subsection{Model}

The problem is formulated in discrete time over $24$ hours using
a constant time step of one hour. The index $t=0,\ldots,T$
(where $T=24$) defines the time discretization grid.

The water volume stored in the dam at time stage $t=0,\ldots,T$, is
a one dimensional random variable $\va{x}_t\in L^2(\omeg,\trib,\prbt;\bbR)$
corresponding to the system state.
This storage variable has to remain between given bounds $\underline{x}$
(minimal volume to be kept in the dam) and $\overline{x}$ (maximal water volume
the dam can contain) so that, for all $t=0,\ldots,T$, the following
almost-sure constraint must hold:
\begin{equation}
\underline{x}\leq \va{x}_t\leq\overline{x} \eqsepv \Pps \eqfinp
\label{dam:bornex}
\end{equation}
The water inflow into the dam at stage $t$ is denoted by $\va{a}_{t}$.
It is a one dimensional random variable with known probability law.
We denote by $\va{u}_t\in L^2(\omeg,\trib,\prbt;\bbR)$ the one dimensional
random variable corresponding to the \emph{desired} volume of water to be
turbinated at stage $t$, and by $\va{e}_{t}$ the \emph{effectively}
turbinated water volume during the same time stage.
In most cases, $\va{e}_t=\va{u}_t$. But this equality is not achievable
if the dam goes under its minimal volume.
Therefore, we have for all $t=0,\ldots,T-1$:
\begin{equation}
\va{e}_t=\min(\va{u}_t,\va{x}_t+\va{a}_{t+1}-\underline{x})
\eqsepv \Pps \eqfinp
\label{dam:commande}
\end{equation}
The shift in the time indices is due to the fact that the decision
to turbinate is supposed to be made before observing the water inflow
volume entering the dam: we are in a \emph{decision-hazard} framework.
Moreover, the control variables $\va{u}_t$ are subject to the following
bounds: for all$t=0,\ldots,T-1$,
\begin{equation}
\underline{u}\leq \va{u}_t\leq\overline{u} \eqsepv \Pps \eqfinp
\label{dam:borneu}
\end{equation}
Taking into account a possible overflow, the dam dynamics write
for $t=0,\ldots,T-1$:
\begin{equation}
\va{x}_{t+1}=\min(\va{x}_t-\va{e}_t+\va{a}_{t+1},\overline{x})
\eqsepv \Pps \eqfinp
\label{dam:dynamique}
\end{equation}
Note that constraints~\eqref{dam:bornex} are fully taken into account
in the modelling of the dam dynamics.
An electricity production $\va{p}_t$ is associated with the effectively
turbinated water volume, and it also depends on the water storage (indeed,
on the water level in the dam, due to the fall height effect):
\begin{equation}
\va{p}_t=g(\va{x}_t,\va{e}_t) \eqfinp
\label{dam:production}
\end{equation}
Let $(\va{d}_{t})_{t=1,\ldots,T}$ denotes the electricity demand,
which is supposed to be a stochastic process with known probability law.
In our decision-hazard framework, production $\va{p}_t$
has to meet demand $\va{d}_{t+1}$: either $\va{p}_t\geq \va{d}_{t+1}$
and the production excess is sold on the electricity market, or
$\va{p}_t\leq \va{d}_{t+1}$ and the gap must be compensated
for either by buying power on the market or by paying a penalty.
The associated cost is modelled~as
\begin{equation}
c_t(\va{d}_{t+1}-\va{p}_t) \eqfinp
\label{dam:cout}
\end{equation}
Ultimately, we suppose that a penalty function $K$ on the final stock
$\va{x}_T$ is given, and that the initial condition $\va{x}_0$ is
a random variable with known probability law.

Let $(\va{w}_t)_{t=0,\ldots,T}$ be the noise random process defined as
\begin{align*}
& \va{w}_0 = \va{x}_0 \eqfinv \\
& \va{w}_t = (\va{a}_t,\va{d}_t) \eqsepv \forall t=1,\ldots,T \eqfinp
\end{align*}
We assume that the noises are fully observed in a non-anticipative way,
and that the control variables are measurable with respect to
the past noises. The dam management problem is then the following.
\begin{subequations}
\begin{equation}
\min_{(\va{u}_{t},\va{x}_{t})} \quad
\bgesp{\sum_{t=0}^{T-1} c_t\big(\va{d}_{t\!+\!1}-g(\va{x}_t,\va{e}_t)\big) +
                        K(\va{x}_T)} \eqfinv
\end{equation}
subject to the constraints
\eqref{dam:commande}--\eqref{dam:borneu}--\eqref{dam:dynamique}
and to the measurability constraints
\begin{equation}
\va{u}_t\preceq (\va{w}_0,\ldots,\va{w}_t) \eqsepv \forall t=0,\ldots,T-1 \eqfinp
\end{equation}
\label{SPDam}
\end{subequations}
It precisely corresponds to the stochastic optimal control problem
formulation~\eqref{SOC} when using the following notations:
\begin{itemize}
\item $w=(a,d)$,
\item $L_t(x,u,w) =
       c_t\big(d-g\big(x,\min(u,x+a-\underline{x})\big)\big)$,
\item $f_t((x,u,w) =
       \min\big(x-\min(u,x+a-\underline{x})+a,\overline{x}\big)$
\item $\Gammaas_t=[\underline{u},\overline{u}]$.
\end{itemize}

\begin{remark}
Both equations~\eqref{dam:commande} and ~\eqref{dam:dynamique} incorporate
the non differentiable operator $\min$. We approximate this non-smooth
operator by the following operator depending on a smoothing parameter $c$:
\begin{equation*}
\min(x,y) \approx
\begin{cases}
y &
    \text{ if } \; y\leq x-c \eqfinv \\
\frac{x+y}{2}-\frac{(x-y)^2}{4c}-\frac{c}{4}&
    \text{ if } \; x-c\leq y\leq x+c \eqfinv \\
x &
    \text{ if } \; y \geq x+c \eqfinp
\end{cases}
\end{equation*}
in order to recover a differentiable problem.
\end{remark}

\subsection{Numerical and functional data}

\begin{figure}[htbp]
\begin{minipage}{6.4cm}
\begin{center}
\includegraphics[width=6.4cm,clip]{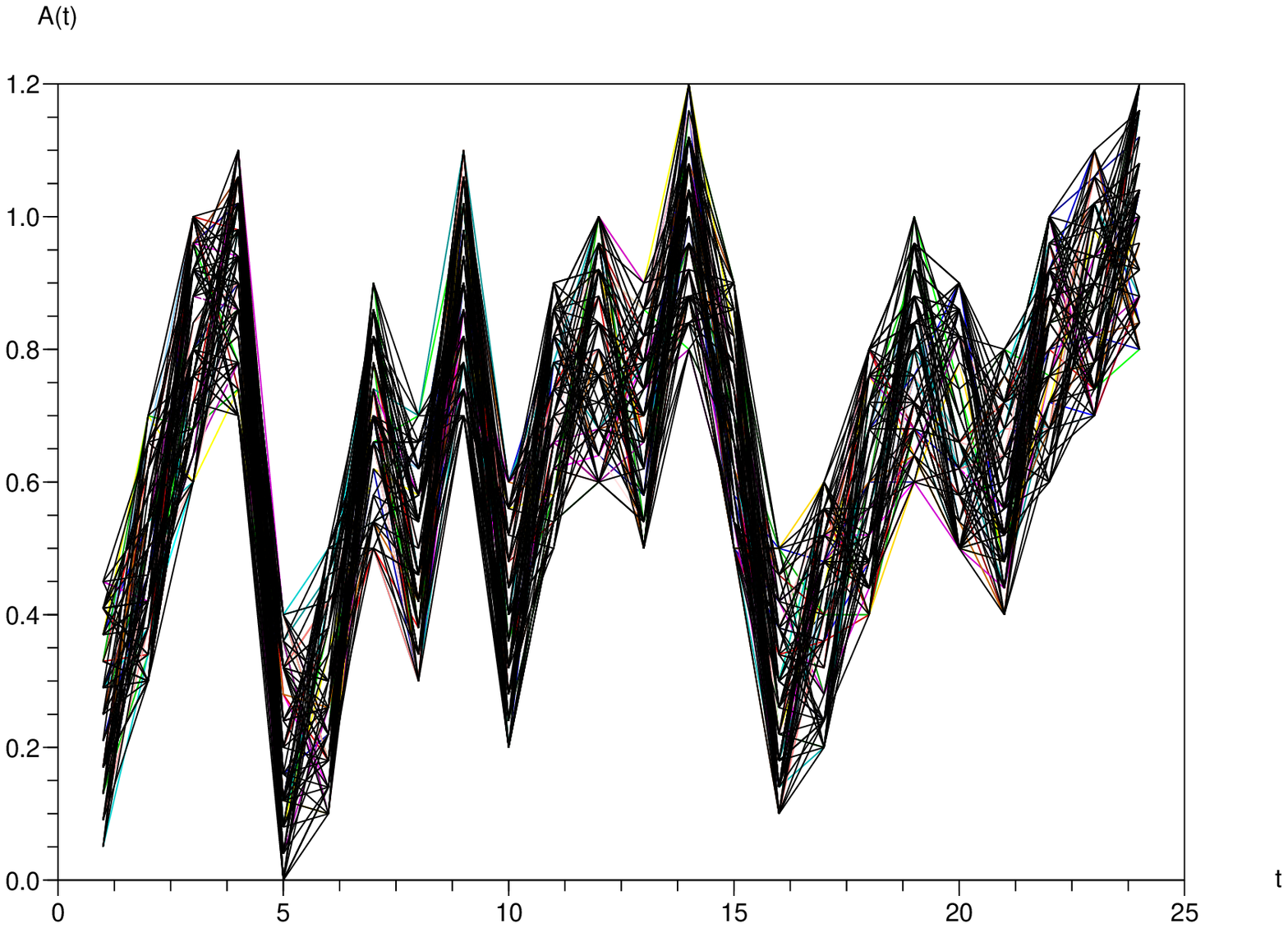}
\caption{\label{dam:apports}Water inflows trajectories}
\end{center}
\end{minipage}
\begin{minipage}{6.4cm}
\begin{center}
\includegraphics[width=6.4cm,clip]{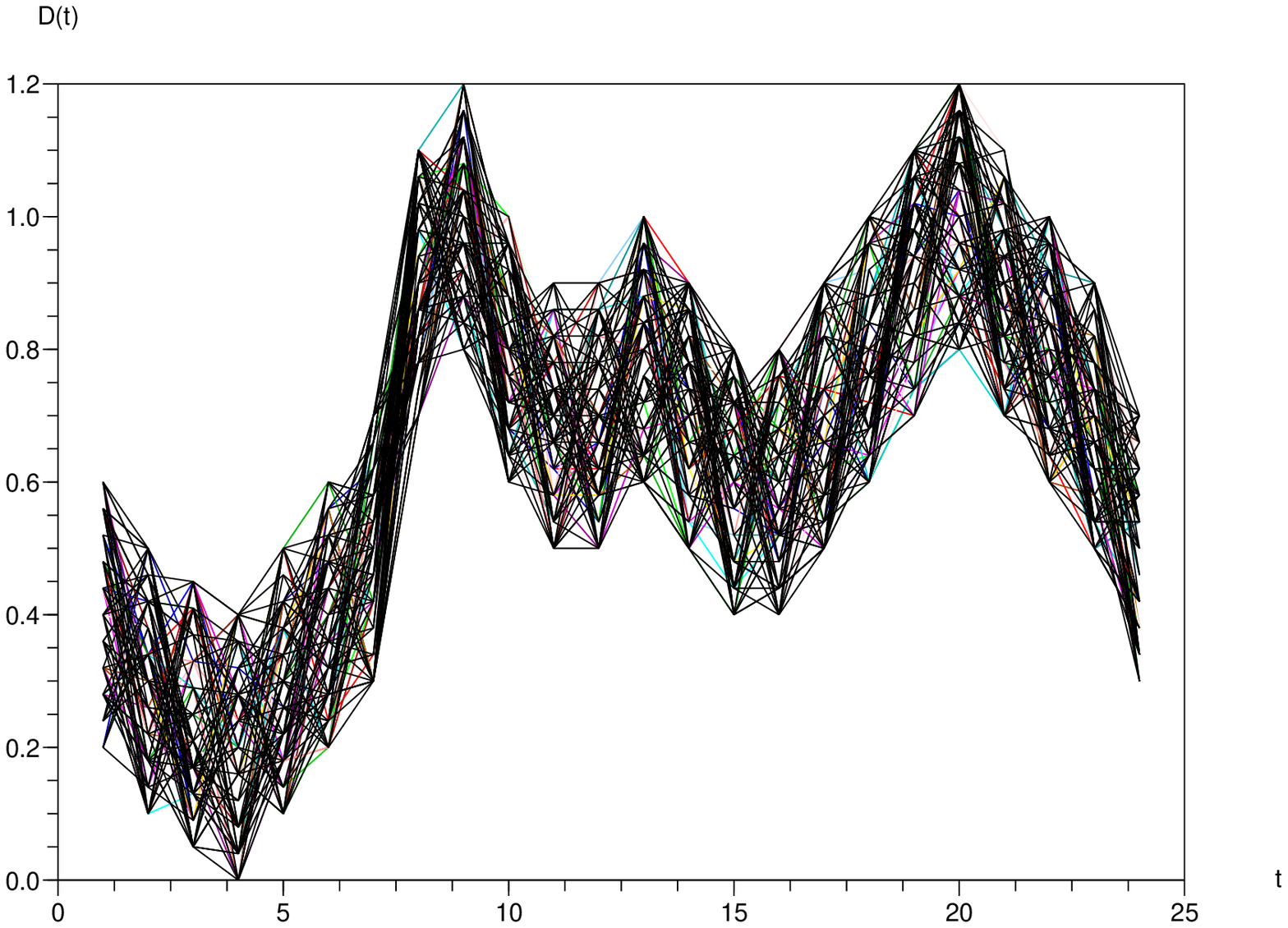}
\caption{\label{dam:demandelec}Electricity demand trajectories}
\end{center}
\end{minipage}
\end{figure}

Both electricity demand and water inflows correspond to white noises,
obtained by adding a discrete disturbance around their mean trajectories.
Using the Monte Carlo method, we draw $N=200$ inflow trajectories
and demand trajectories, which are depicted in Figure~\ref{dam:apports}
and Figure~\ref{dam:demandelec} respectively, the associated particles
being denoted
$\big(\va{a}_t^k\big)_{t=1,\ldots,T}^{k=1,\ldots,N}$ and
$\big(\va{d}_t^k\big)_{t=1,\ldots,T}^{k=1,\ldots,N}$.

The initial state $\va{x}_0$ follows a uniform probability law
over $[\underline{x},\overline{x}] = [0,2]$. We also draw
$N$ particles $(\va{x}_0^k)^{k=1,\ldots,N}$ for the initial state
and each one is associated with the previous trajectories with
the same index~$k$ to form one scenario among~$N$.
The control random variables $\va{u}_t$ are subject for each
$t=0,\ldots,T-1$ to the bounds $[\underline{u},\overline{u}] = [0,1]$.

The mapping $g$ modelling the electricity production $\va{p}_t$
is chosen to represent a linear variation  between $0.5$ and $1$
with respect to the water fall height $\va{x}_t-\underline{x}$:
\begin{equation*}
g(\va{x}_t,\va{e}_t) = \va{e}_t
\frac{\va{x}_t+\overline{x}-2\underline{x}}{2(\overline{x}-\underline{x})}
\eqfinp
\end{equation*}

\begin{figure}[htbp]
\begin{center}
\includegraphics[width=7cm,clip]{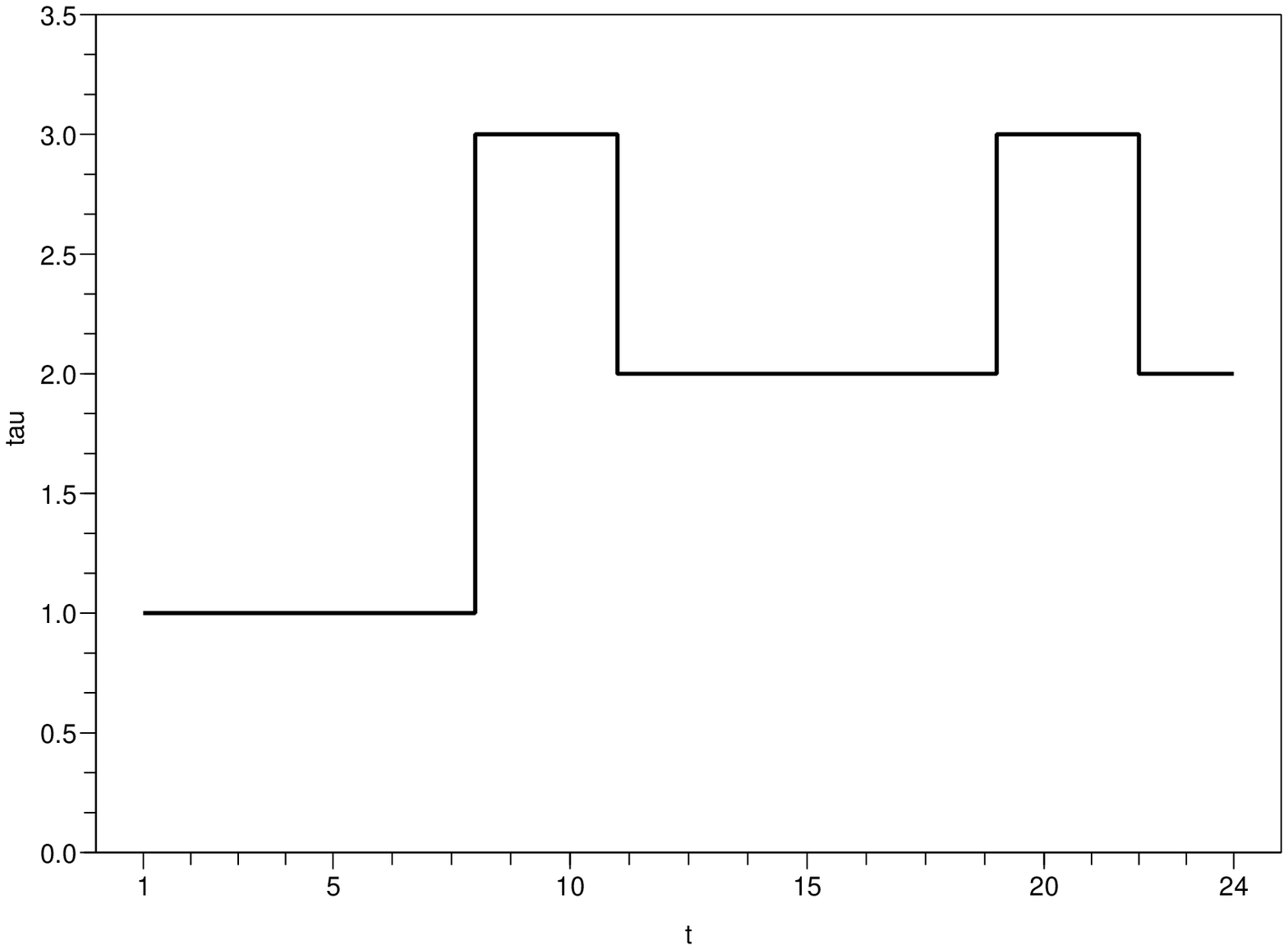}
\caption{\label{dam:prixelec}Electricity price $\tau_t$}
\end{center}
\end{figure}

The expression of the instantaneous cost $c_{t}$ is
\begin{equation*}
c_t(y)=\tau_{t}(e^y-1) \eqfinv
\end{equation*}
where $\tau_{t}$ is the electricity price at stage $t$.
The variation of this price is depicted in Figure~\ref{dam:prixelec}.

The final cost is an incentive to fill the dam at the end of the day:
\begin{equation*}
K(x)=12(x-\overline{x})^2.
\end{equation*}

\subsection{Results using Dynamic Programming}

We apply Algorithm~\ref{DPApprox} to Problem~\eqref{SPDam} using
an even discretization of the state space $[\underline{x},\overline{x}]$
in $n=200$ points:
\begin{equation*}
x_i=\underline{x}+\frac{i-1}{n-1}(\overline{x}-\underline{x})
\eqsepv \forall i=1,\ldots,n \eqfinp
\end{equation*}
The two operators $\gridop{R}_\bbR$ and $\gridop{R}_{\espacea{U}_t}$
are linear interpolation operators: to compute the value of a function
outside the grid, we consider the weighted mean of the two surrounding
grid points. The optimal feedback laws $\phi_t$ obtained at each time
stage $t=0,\ldots,T-1$ will be used as a reference in the comparison
with the other resolution methods. The optimal cost is obtained by
simulating the system using the optimal feedback laws over
all trajectories:
\begin{equation}
c \defegal \besp{V_0(\va{x}_0)} =
           \frac{1}{2}\int_0^2 V_0(x)\mathrm{d}x =
           6.48 \eqfinp
\label{cSDP}
\end{equation}

\subsection{Results obtained by Stochastic Programming}

We then make use of a stochastic programming technique to solve
Problem~\eqref{SPDam}. Using quantization techniques, we first
generate a scenario tree from the $N=200$ noises trajectories.
We will not discuss here the quantization method used to build such
a scenarios tree and refer to \cite{TheseBarty} for further details.
The resulting tree includes $2$ nodes at stage $t=0$, $4$ nodes
at stage $t=1$ and so on till stage $t=6$ for which we have
$2^{6+1}=128$ nodes. As $2^{8}>200$, the tree structure becomes
deterministic as soon as $t\geq 7$, each node in the tree corresponding
to $t \geq 7$ having a unique future (Figure~\ref{treegen} just sketches
the beginning of the story).

\begin{figure}[htbp]
\begin{minipage}{6.45cm}
\begin{center}
\includegraphics[width=6cm,clip]{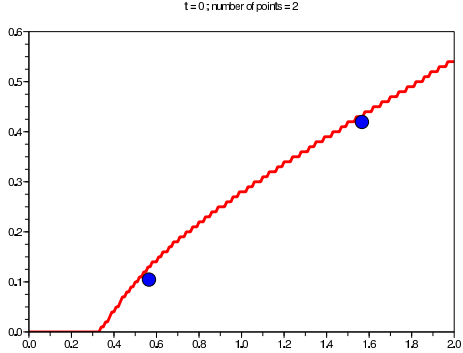}
\caption{\label{feed0Tree}Scenario tree: optimal control ($t=0$)}
\end{center}
\end{minipage}
\begin{minipage}{6.45cm}
\begin{center}
\includegraphics[width=6cm,clip]{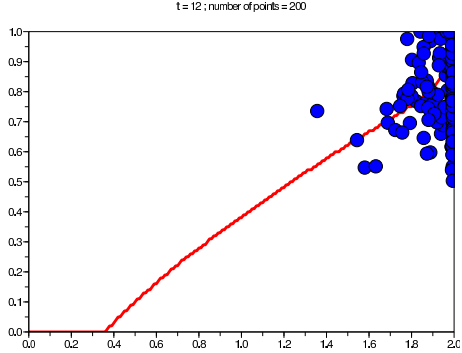}
\caption{\label{feed12Tree}Scenario tree: optimal control ($t=12$)}
\end{center}
\end{minipage}
\end{figure}

\begin{figure}[htbp]
\begin{center}
\includegraphics[width=6cm,clip]{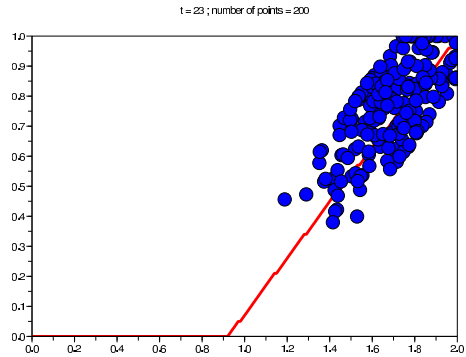}
\caption{\label{feed23Tree}Scenario tree: optimal control ($t=23$)}
\end{center}
\end{figure}

Problem~\eqref{SPDam} is then formulated and optimized over the tree:
the optimization process yields a pair $(x_{\nu},u_{\nu})$ of optimal
values for the state and the control at each node $\nu$ of the tree.
The next figures depict the optimal pairs of particles at different
time stages (represented by dots) and the optimal feedback laws obtained
by Dynamic Programming (represented by continuous curves).
The comparison leads to the following conclusions.
\begin{itemize}
\item There are only two nodes corresponding to $t=0$ in the scenario tree,
and therefore only two optimal control particles. These two particles fit
the optimal feedback obtained by Dynamic Programming rather accurately
(see  Figure~\ref{feed0Tree}), but it would be difficult to synthesize
a feedback law with such a limited number of points.
\item At stages $t=12$ and $t=23$, $200$ optimal control particules
are available. Nevertheless these particles have a visible huge variance
(see Figures~\ref{feed12Tree} and \ref{feed23Tree}), so that it would again
be difficult to synthesize a feedback law.
\end{itemize}

\subsection{Results with the adaptive mesh method}

\begin{figure}[htbp]
\begin{minipage}{6.45cm}
\begin{center}
\includegraphics[width=6.45cm,clip]{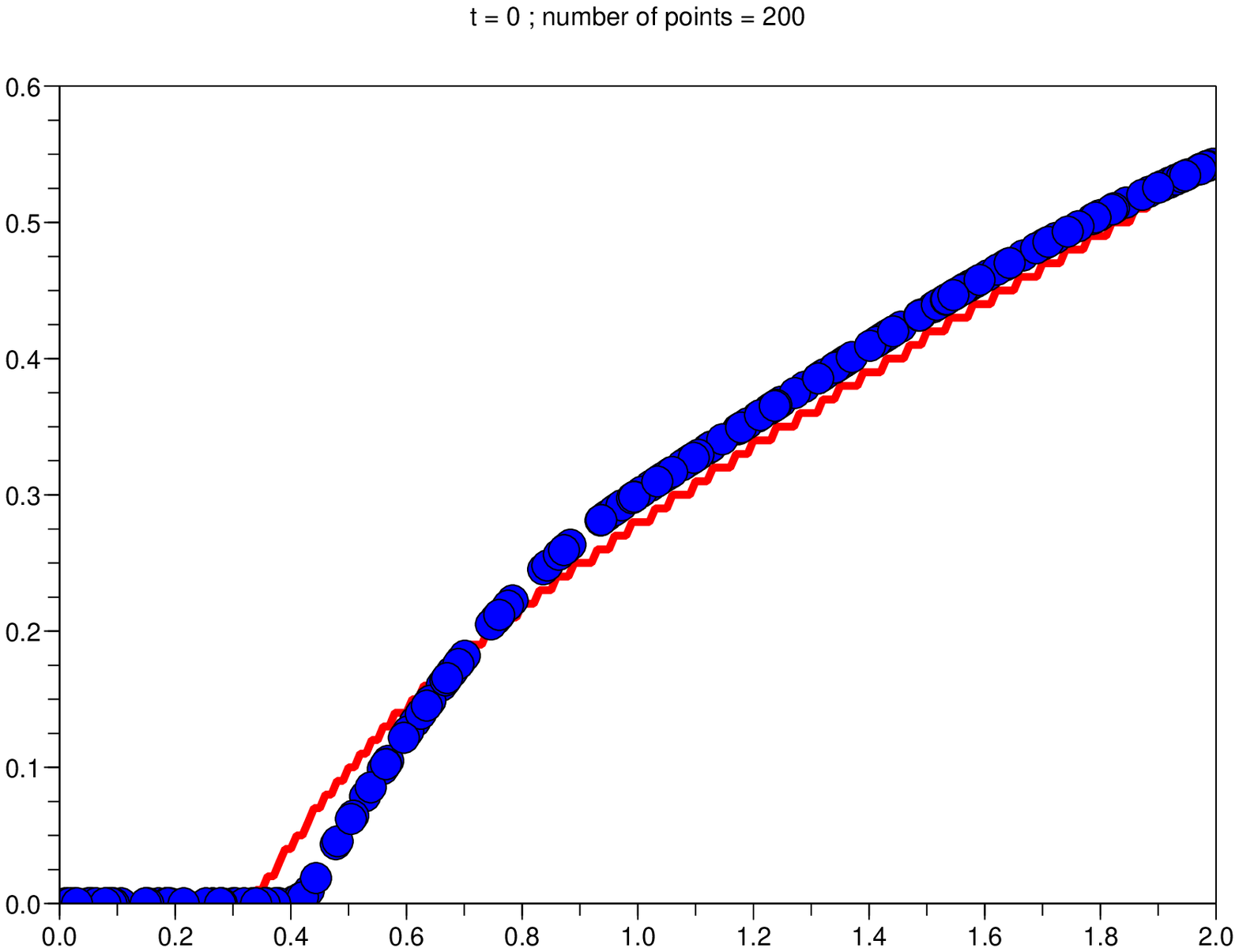}
\caption{\label{feed0Part}Particle: optimal control ($t=0$)}
\end{center}
\end{minipage}
\begin{minipage}{6.45cm}
\begin{center}
\includegraphics[width=6.45cm,clip]{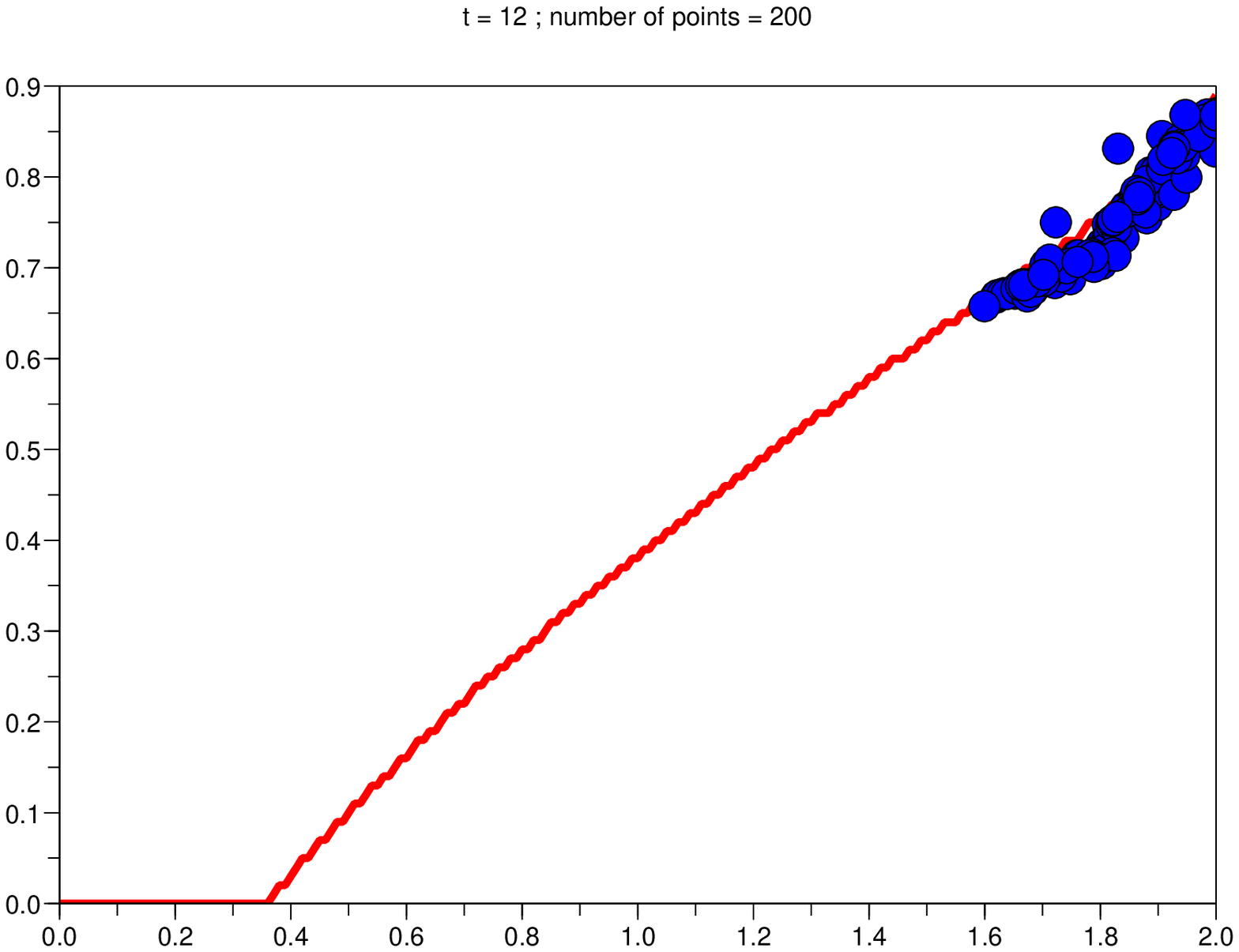}
\caption{\label{feed12Part}Particle: optimal control ($t=12$)}
\end{center}
\end{minipage}
\end{figure}

\begin{figure}[htbp]
\begin{center}
\includegraphics[width=6.45cm,clip]{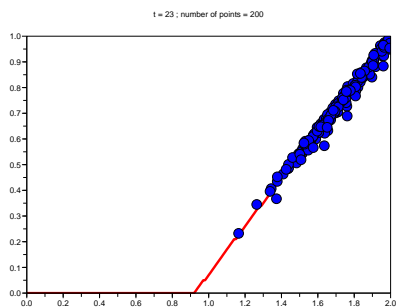}
\caption{\label{feed23Part}Particle: optimal control ($t=23$)}
\end{center}
\end{figure}

We ultimately apply Algorithm~\ref{gradproj} to solve the
hydro-electric dam management problem~\eqref{SPDam}.
The result of this algorithm is an optimal feedback
law $\phi_t$ for every time stage $t=0,\ldots,T-1$.
We then draw new noise trajectories independent from those
used by the algorithm, and we simulate the system behavior
along these new trajectories using the optimal feedback laws:
\begin{equation*}
\begin{split}
&\va{x}_0=\va{w}_0 \eqfinv \\
%& \va{u}_t=\phi_t(\va{x}_{t})
%  \eqsepv \forall t=0,\ldots,T-1 \eqfinv \\
& \va{x}_{t+1}=
  f_t\big(\va{x}_{t},
          \phi_t(\va{x}_{t}),
          \va{w}_{t+1}\big)
  \eqsepv \forall t=0,\ldots,T-1 \eqfinv
\end{split}
\end{equation*}
and thus obtain an approximation of the optimal cost generated
by this algorithm:
\begin{equation*}
c = 6.51 \approx \bgesp{\sum_{t=0}^{T-1}
   L_t\big(\va{x}_t,
           \phi_t(\va{x}_{t}),
           \va{w}_{t+1}\big) +
   K(\va{x}_T)} \eqfinp
\end{equation*}
This optimal cost is close to the cost generated by Dynamic Programming.
We are also interested in the controls generated by the adaptive mesh
method. To this purpose, Figures~\ref{feed0Part}, \ref{feed12Part} and
\ref{feed23Part} show the optimal control particles given by the adaptive
mesh method (dots), to be compared with the optimal feedback laws
obtained by Dynamic Programming (curves) for different time stages.

We note that the optimal particles obtained by the adaptive mesh
method are close to the feedback laws obtained by Dynamic Programming.
By construction, there are the same number of particles at each
time stage $t$, and the dispersion of the particles remains at first
sight constant from the beginning to the end of the time horizon.
This represents a significant advance compared to the scenario
tree method. On the other hand, observe that particles may sometimes
concentrate in restricted parts of the state space
(see Figure~\ref{feed12Part}): in our view, this is not a \emph{drawback},
but an \emph{advantage} of the proposed method in that the optimal
feedback is computed only where it is needed to do so. Indeed, the
particles distribute adaptively and automatically according to
the optimal probability density of the state (we call this an
``adaptive mesh'' --- in Figure~\ref{feed0Part}, the distribution
is even because the initial condition is uniformly distributed),
and this is an advantage over Dynamic Programming in which a uniform
grid is defined a priori over the whole state space, irrespective of
the optimal solution distribution.

\section{Conclusions and perspectives}

In this paper we presented new tractable methods for solving
stochastic optimal control problems in the discrete time case.
We derived several forms of the optimality conditions for such problems:
\begin{itemize}
\item non-adapted optimality conditions~\eqref{eqPMS1} as well as adapted
optimality conditions~\eqref{eqPMS2} with measurability constraints
on the past noise variables; these conditions incorporate conditional
expectations and the dimension of the conditioning random variable
grows with the number of time stages;
\item in the Markovian case, non-adapted optimality conditions~\eqref{eqPMS1X}
and adapted optimality conditions~\eqref{eqPMS2X} with measurability
constraints on the state variable; the conditional expectations are
taken with respect to the instantaneous state variable whose dimension
is constant over the time stages, but the mean square error when
approximating such conditional expectations depends on the state
space dimension;
\item still in the Markovian case, functional optimality
conditions~\eqref{eqPMS3X} including only expectations.
\end{itemize}
The last conditions have been used to devise a gradient-like adaptive
mesh algorithm in order to solve stochastic optimal control problems
in the Markovian case, and we applied the algorithm to a hydro-electric
dam management problem.

In light of the numerical results, it is clear that the proposed
adaptive mesh algorithm represents a significant advance with respect
to usual stochastic programming techniques (same number of particles
and same particle dispersion at every time stage). In addition, the
adaptive mesh feature may save useless computations in some problems,
depending on the profile of the optimal state probability density.
In fact, the only a priori discretization concerns noise particles,
which does not depend on the dimension of the underlining state space:
the only operator which could be dimension-dependent is the interpolation
operator.

Future work will concentrate on the convergence rate of the mesh
algorithm with respect to the number $N$ of noise trajectories.
We will also deal with stochastic optimal control problems involving
a multi-dimensional state vector, and try to quantify the impact
of the interpolation operator on the approximation error.

\appendix
\section{\label{app-OC-Hilbert}Optimization on an Hilbert space: a special case}
Let $\espacef{H}$ be an Hilbert space, let $\Hfe$ be a closed convex
subset of $\espacef{H}$ and let $f$ be a real valued function defined
on $\espacef{H}$. We consider the following optimization problem:
\begin{equation}
\min_{x\in\Hfe} \; f(x) \eqfinp
\label{GOP}
\end{equation}
In the following, $\fcara{H}$ will denote the indicator function
of a subset $H\subset\espacef{H}$, namely
\begin{equation*}
\fcara{H}(x) = \left\{
\begin{array}{ll}
0       & \text{if $x\in H$} \eqfinv\\
+\infty & \text{otherwise} \eqfinp
\end{array}
\right.
\end{equation*}
The optimization literature gives different expressions
for the necessary optimality conditions of an optimization problem
in a general Hilbert space (see e.g. \cite{EkelandTemam76}).
For instance, if $x^\sharp\in\espacef{H}$ is solution of~\eqref{GOP},
then the following statements are equivalent:
\begin{subequations}
\begin{align}
& \forall x\in\Hfe \eqsepv \proscal{f'(x^\sharp)}{x-x^\sharp}\geq 0
  \eqfinv \label{cond1} \\
& f'(x^\sharp)\in -\partial\fcara{\Hfe}(x^\sharp)
  \eqfinv \label{cond2} \\
& \forall \epsilon> 0 \eqsepv
  x^\sharp=\proj{\Hfe}{x^\sharp-\epsilon \gradi{f}(x^\sharp)}
  \eqfinp \label{cond3}
\end{align}
\label{cond123}
\end{subequations}
We now consider a specific structure for the feasible set $\Hfe$.
More precisely, we assume that $\Hfe=\Hcv\cap\Hsp$, $\Hsp$ being
a \emph{closed subspace} of $\espacef{H}$ and $\Hcv$ being
a \emph{closed convex subset} of $\espacef{H}$. We moreover
assume that the following property holds.

\vspace{0.2cm}

\begin{assumption}
The sets $\Hsp$ and $\Hcv$ are such that
$\proj{\Hcv}{\Hsp}\subset\Hsp$.
\label{ass:proj-inter}
\end{assumption}

\vspace{0.2cm}

\noindent
Then the projection operator on~$\Hfe$ has the following property.

\begin{lemma}
Under Assumption \ref{ass:proj-inter}, the following
relation holds true:
\begin{equation*}
\projop{\Hcv\cap\Hsp}=\projop{\Hcv}\circ\projop{\Hsp} \eqfinp
\end{equation*}
\label{lem:proj-inter}
\end{lemma}

\begin{proof}
Let $y\in\Hcv\cap\Hsp$. Then
\begin{multline*}
\proscal{x-\proj{\Hcv}{\proj{\Hsp}{x}}}
        {y- \proj{\Hcv}{\proj{\Hsp}{x}}} = \\
\proscal{x-\proj{\Hsp}{x}}
        {y- \proj{\Hcv}{\proj{\Hsp}{x}}} + \\
\proscal{\proj{\Hsp}{x}-\proj{\Hcv}{\proj{\Hsp}{x}}}
        { y- \proj{\Hcv}{\proj{\Hsp}{x}}} \eqfinp
\end{multline*}
From the characterization of the projection
of $z\defegal \proj{\Hsp}{x}$ over the convex subset $\Hcv$,
the last inner product in the previous expression is non
positive: therefore,
\begin{equation*}
\proscal{x-\proj{\Hcv}{z}}{y- \proj{\Hcv}{z}}
\leq \proscal{x- z}{y- \proj{\Hcv}{z}} \eqfinp
\end{equation*}
From $\proj{\Hcv}{\Hsp}\subset\Hsp$, we deduce that
$y-\proj{\Hcv}{z}\in\Hsp$. Since $\projop{\Hsp}$
is a self-adjoint operator, we have
\begin{align*}
\proscal{x-\proj{\Hcv}{z}}{y-\proj{\Hcv}{z}}
  &\leq \proscal{x- z}{\proj{\Hsp}{y-\proj{\Hcv}{z}}} \eqfinv \\
  &\leq \proscal{\proj{\Hsp}{x-z}}{y-\proj{\Hcv}{z}} \eqfinv \\
  &\leq 0 \eqfinv
\end{align*}
the last inequality arising from the fact that
$\proj{\Hsp}{x-z} = \proj{\Hsp}{x}-\proj{\Hsp}{z}$
$=0$ (since $\Hsp$ is a linear subspace, then
$\proj{\Hsp}{\cdot}$ is a linear operator). We thus conclude
that, for all $y\in\Hcv\cap\Hsp$,
\begin{equation*}
\proscal{x-\projop{\Hcv}\circ\proj{\Hsp}{x}}
        {y- \projop{\Hcv}\circ\proj{\Hsp}{x}}\leq 0 \eqfinv
\end{equation*}
a variational inequality which characterizes
$\projop{\Hcv}\circ\proj{\Hsp}{x}$
as the projection of $x$ over $\Hfe=\Hcv\cap\Hsp$. \qquad
\end{proof}

The following proposition gives necessary optimality conditions
for Problem~\eqref{GOP} when the feasible set $\Hfe$
has the specific structure $\Hcv\cap\Hsp$.

\begin{proposition}
We suppose that Assumption \ref{ass:proj-inter} is fulfilled
and that $f$ is differentiable.
If $x^\sharp$ is solution of~\eqref{GOP}, then
\begin{equation*}
\proj{\Hsp}{f'(x^\sharp)}\in -\partial\fcara{\Hcv}(x^\sharp) \eqfinp
\end{equation*}
\label{propCOpr}
\end{proposition}

\begin{proof}
Let $x^\sharp$ be solution of~\eqref{GOP}. Using Condition~\eqref{cond3}
and Lemma \ref{lem:proj-inter}, we obtain that
\begin{equation*}
x^\sharp=\projop{\Hcv}\circ\proj{\Hsp}{x^\sharp-\epsilon \gradi{f}(x^\sharp)}
\eqsepv \forall \epsilon \geq 0 \eqfinp
\end{equation*}
But $\projop{\Hsp}$ is a linear operator and $x^\sharp\in\Hsp$, so that
\begin{equation*}
x^\sharp=\proj{\Hcv}{x^\sharp-\epsilon\:\proj{\Hsp}{\gradi{f}(x^\sharp)}}
\eqfinp
\end{equation*}
From~\eqref{cond123}, the last relation is equivalent to
$\proj{\Hsp}{f'(x^\sharp)}\in -\partial\fcara{\Hcv}(x^\sharp)$. \qquad
\end{proof}

\bibliographystyle{siam}

\bibliography{biblio}

\begin{thebibliography}{10}

\bibitem{AubinFrankowska90}
{\sc J.-P. Aubin and H.~Frankowska}, {\em Set-valued Analysis}, Birkh\"auser,
  Boston, 1990.

\bibitem{TheseBarty}
{\sc K.~Barty}, {\em Contributions \`{a} la discr\'{e}tisation des contraintes
  de mesurabilit\'{e} pour les probl\`{e}mes d'optimisation stochastique},
  {PhD} dissertation, \'{E}cole Nationale des Ponts et Chauss\'{e}es, 2004.

\bibitem{barty03}
{\sc K.~Barty, P.~Carpentier, J.-P. Chancelier, G.~Cohen, M.~De~{L}ara, and
  T.~Guilbaud}, {\em Dual effect free stochastic controls}, Annals of
  Operations Research, 142 (2006), pp.~41--62.

\bibitem{Bellman:1957}
{\sc R.~Bellman}, {\em Dynamic Programming}, Princeton University Press, New
  Jersey, 1957.

\bibitem{Bertsekas:1976}
{\sc D.~Bertsekas}, {\em Dynamic programming and stochastic control}, Acad.
  Press, 1976.

\bibitem{Bertsekas-Shreve:1996}
{\sc D.~Bertsekas and S.~Shreve}, {\em Stochastic Optimal Control: the
  discrete-time case}, Athena Scientific, Belmont, 1996.

\bibitem{Brei92}
{\sc L.~Breiman}, {\em Probability}, Society for Industrial and Applied
  Mathematics, Philadelphia PA, 1992.

\bibitem{BrodieGlasserman97}
{\sc P.~Brodie, M.~Glasserman}, {\em A stochastic mesh method for pricing high
  dimensional american options}, The journal of computational finance, 7
  (2004).

\bibitem{TheseDallagi}
{\sc A.~Dallagi}, {\em M\'{e}thodes particulaires en commande optimale
  stochastique}, {PhD} dissertation, {U}niversit\'{e} {P}aris {I}
  {P}anth\'{e}on-{S}orbonne, 2007.

\bibitem{DupacovaGroweKuskaRomisch03}
{\sc J.~Dupa\`cov\'a, N.~Gr\"owe-Kuska, and W.~R\"omisch}, {\em Scenario
  reduction in stochastic programming. {A}n approach using probability
  metrics}, Math. Program., 95 (2003), pp.~493--511.

\bibitem{EkelandTemam76}
{\sc I.~Ekeland and R.~Temam}, {\em Convex analysis and variational problems},
  SIAM, Philadelphia, 1999.

\bibitem{HeitschRomisch03}
{\sc H.~Heitsch and W.~R\"omisch}, {\em Scenario reduction algorithms in
  stochastic programming}, Comput. Optim. Appl.,  (2003), pp.~187--206.

\bibitem{Hiriart82}
{\sc J.-B. Hiriart-Urruty}, {\em Extension of lipschitz integrands and
  minimization of nonconvex integral functionals : {A}pplications to the
  optimal recourse problem in dicrete time}, Probability and mathematical
  statistics, 3 (1982), pp.~19--36.

\bibitem{Leese:1974}
{\sc S.~Leese}, {\em Multifunctions of {S}ouslin type}, Bull. Austral. Math.
  Soc., 11 (1974), pp.~395--411.

\bibitem{Romisch05}
{\sc J.~Outrata and W.~R\"omisch}, {\em On optimality conditions for some
  nonsmooth optimization problems over {$L^p$} spaces}, Journal of Optimization
  Theory and Applications, 126 (2005), pp.~411--438.

\bibitem{Pflug01}
{\sc G.~Pflug}, {\em Scenario tree generation for multiperiod financial
  optimization by optimal discretization}, Math. Program., 89 (2001),
  pp.~251--271.

\bibitem{Rao04}
{\sc M.~Rao}, {\em Measure Theory and Integration}, {P}ure and {A}pplied
  {M}athematics {S}eries, {M}arcel {D}ekker {I}nc, 2004.

\bibitem{Rockafellar-Wets:1998}
{\sc R.~Rockafellar and R.-B. Wets}, {\em Variational Analysis}, Springer
  Verlag, Berlin Heidelberg, 1998.

\bibitem{TheseStrugarek}
{\sc C.~Strugarek}, {\em Approches variationnelles et autres contributions en
  optimisation stochastique}, {PhD} dissertation, \'{E}cole Nationale des Ponts
  et Chauss\'{e}es, 2006.

\bibitem{ThenieVial06}
{\sc J.~Th\'{e}ni\'{e} and J.~Vial}, {\em Step decision rules for multistage
  stochastic programming: a heuristic approach}.
\newblock \href{http://www.optimization-online.org/DB_HTML/2006/08/1440.html}
  {http://www.optimization-online.org/DB\_HTML/2006/08/1440.html}, 2006.

\bibitem{Wagner:1977}
{\sc D.~Wagner}, {\em Survey of measurable selection theorems}, SIAM J. Control
  Optim., 15 (1977), pp.~859--903.

\end{thebibliography}

\end{document}